\documentclass[10pt]{article}

\usepackage{amsmath,amssymb,amsfonts}

\usepackage{tikz}
\usepackage{listings}
\usepackage{subfigure}
\usepackage{curves}
\usepackage{xypic}  
\usepackage[mathscr]{eucal}
\usepackage[pdftex]{hyperref}
\usepackage[all]{xy}
\usepackage{graphics,graphicx}
\usepackage{epsfig}
\usepackage{color}
\definecolor{red-}{rgb}{1.0,0.0,0.0}
\definecolor{grey}{rgb}{0.6, 0.6, 0.6}
\definecolor{brown}{rgb}{0.5,0.2,0.0}
\definecolor{brown-}{rgb}{0.0,0.1,1.0}
\definecolor{green-}{rgb}{0.0, 0.6, 0.0}
\definecolor{gold}{rgb}{0.8,0.7,0.0}
\definecolor{black}{rgb}{0.0,0.0,0.0}
\definecolor{DarkGreen}{rgb}{0.0,0.3,0.2}
\definecolor{LightGreen}{rgb}{0.8,1.0, 0.8}
\definecolor{yellow}{rgb}{0.9,0.9,0.0}
\definecolor{blue-}{rgb}{0.0,0.1,1.0}

\usepackage{bm}  
\usepackage{cancel}

\font\tenmsb=msbm10
\font\sevenmsb=msbm7
\font\fivemsb=msbm5

\newfam\msbfam
\textfont\msbfam=\tenmsb
\scriptfont\msbfam=\sevenmsb
\scriptscriptfont\msbfam=\fivemsb

\font\teneufm=eufm10
\font\seveneufm=eufm7
\font\fiveeufm=eufm5
\newfam\eufmfam
\textfont\eufmfam=\teneufm
\scriptfont\eufmfam=\seveneufm
\scriptscriptfont\eufmfam=\fiveeufm

\usepackage[all]{xy}

\newcommand\sE{{\cal E}}

\newcommand{\x}{{\bm x}}

\newcommand\R{{\mathbb R}}

\newcommand{\be}{\begin{equation}}
\newcommand{\ee}{\end{equation}}

\newcommand{\epsbold}{{\bm \varepsilon}}
\def\Jac{\mathop{\rm Jac}\nolimits}

\newtheorem{theorem}{Theorem}[section]

\newtheorem{question}[theorem]{Question}

\newtheorem{definition}[theorem]{Definition}
\newtheorem{re}[theorem]{Remark}
\newtheorem{defre}[theorem]{Definition--Remark}
\newtheorem{pargrph}[theorem]{}
\newtheorem{examp}[theorem]{Example}

\newtheorem{MM}[theorem]{ }

\newtheorem{probdef}[theorem]{Problem--Definition}
\newtheorem{ex}[theorem]{Esercizio}
\newtheorem{propdef}[theorem]{Proposition--Definition}

\makeatletter
\ifnum\@ptsize=0\fi
\ifnum\@ptsize=2\fi

\def\cocoa{{\hbox{\rm C\kern-.13em o\kern-.07em C\kern-.13em o\kern-.15em A}}}

\newenvironment{rem*}{\begin{re}\em}{\end{re}}
\newenvironment{example*}{\begin{examp}\em}{\end{examp}}
\newenvironment{definition*}{\begin{definition}\em}{\end{definition}}
\newenvironment{probdef*}{\begin{probdef}\em}{\end{probdef}}
\newenvironment{question*}{\begin{question}\em}{\end{question}}

\newenvironment{prgrph*}[1]{\indent\begin{pargrph}{\bf #1.}\em\
}{\end{pargrph}}
\newenvironment{defre*}{\begin{defre}\em}{\end{defre}}
\newenvironment{MM*}{\begin{MM}\em}{\end{MM}}
\newenvironment{ex*}{\begin{ex}\em}{\end{ex}}
\newenvironment{propdef*}{\begin{propdef}\em}{\end{propdef}}

\begin{document}

\title{Corrigendum for   ``Almost vanishing polynomials \\and an  application  to the Hough transform"
\footnote{2010
{\em Mathematics Subject Classification}. Primary  26C10, 15A60, 14Q10;
Secondary 14H50  \newline
\indent{{\em Keywords and phrases.}  Weighted matrix norm; hypersurface;  crossing area conditions} }}
\author{ Maria-Laura Torrente and Mauro C. Beltrametti }

\date{}
\maketitle


In this note we correct a technical error occurred in \cite{TB}.
 This affects the bounds given in  that paper, 
even though  the structure and the logic of all proofs remain  fully unchanged. The error is due to a repeated wrong use  of 
H\"older's inequality (a transpose of a matrix  was missed). The first time it occurs is  in the first inequality of formula $(8)$ of \cite{TB}.   Indeed, the correct version of that formula is:

\begin{eqnarray*}\label{diseq2}
 |\Jac_f(p) (p^*-p)^t| & = &  |\Jac_f(p) \sE^{-1} \sE (p^*-p)^t| \nonumber \\ 
 &\le&  \|\big(\Jac_f(p) \sE^{-1}\big)^t\|_1 \|\sE (p^*-p)^t\|_\infty \nonumber \\ 
 &\le& \|\sE^{-1}\Jac_f(p)^t \|_1 \le  \|\Jac_f(p)^t\|_1 \|\sE^{-1}\|_1.
\end{eqnarray*}


 We refer to our paper \cite{TB}, using the same notation.
Here, we confine ourselves to state the correct versions of 
 Proposition 2.1,  Proposition 2.5, Proposition 3.2, and Proposition 4.3 of \cite{TB}, respectively. Up to the error pointed out and corrected as above, the proofs go parallel to those in~\cite{TB}. In the following, $P$  denotes the multivariate polynomial ring $\R[\x]=\R[x_1,\ldots,x_n]$.

\bigskip
\bigskip
\bigskip
\noindent{\bf Proposition 2.1}
{\em Let $f=f(\x)$ be a non-constant polynomial of~$P$,
let $p$ be a point of~$\R^n$, and let ${\bf C}(p) \subseteq {\bf B}(p)$
be an $(\infty,\epsbold)$-unit cell   centered at $p$. If
$$
| f(p)| > \|\Jac_f(p)^t\|_1 \epsbold_{\rm max} +\frac{n}{2}{\mathsf H} \epsbold_{\rm max}^2 =:B_1,
$$
then the hypersurface  of equation $f=0$ does not cross ${\bf C}(p).$
}
\bigskip

\noindent{\bf Proposition 2.5}
{\em Let $f(\x)$ be a degree $\geq 2$   polynomial of~$P$.
Let $p$ be a point of $\R^n$ and let
${\bf C}(p) \subseteq {\bf B}(p)$
be an $(\infty,\epsbold)$-unit cell  centered at $p$. If
$$
|f(p)| > \|{\rm Jac}_f(p)^t\|_1 \epsbold_{\rm max} + \frac{n}{2}  \|H_f(p)\|_\infty \epsbold_{\rm max}^2:=B_1',
$$
then  the hypersurface of equation $f=0$ does not cross ${\bf C}(p)$ neglecting contributions of order ${\rm O}(\epsbold_{\rm max} ^3)$.
}
\bigskip

\noindent{\bf Proposition 3.2}
{\em Let $f=f(\x)$ be a degree $\geq 2$ polynomial of~$P$,
let $p$ be a point of $\R^n$ such that $\Jac_f(p)$ is not the zero vector, 
and let ${\bf C}(p) \subseteq {\bf B}(p)$
be an $(\infty,\epsbold)$-unit cell  centered at $p$.
Let $R$ be a positive real number such that 
$R < \min \big\{\epsbold_{\rm min}, \frac{\|\Jac_f(p)\|_1}{{\mathsf H}}\big\}$.
Set $c:=\max\{2, \sqrt{n}\}$. If 
$$
|f(p)| < \frac{2R}{\mathsf J(c+ n^{5/2}{\mathsf H} {\mathsf J} R)}  =: B_2,
$$
then the hypersurface of equation $f=0$ crosses ${\bf C}(p)$.
}

\bigskip
\noindent{\bf Proposition 4.3}
{\em Let $f=f(\x)$ be a degree $\geq 2$ polynomial of~$P$,
let $p$ be a point of $\R^n$ such that  Jacobian $\Jac_f(p)$ and the Hessian matrix $H_f(p)$ are nontrivial, 
and let ${\bf C}(p) \subseteq {\bf B}(p)$
be an $(\infty,\epsbold)$-unit cell  centered at $p$.
Let $R$ be a positive real number such that
$R < \min \left \{\epsbold_{\rm min}, \frac{\|\Jac_f(p)\|_1}{n^2 \|H_f(p)\|_\infty} \right\}$, 
let $c:=\max\{2, \sqrt{n}\}$ and set 
$$
\Theta:=  \|\Jac_f^\dagger(p)\|_\infty+  n^2(1+2 \sqrt{n}) \frac{\|H_f(p)\|_\infty}{\|\Jac_f(p)\|_1^2} R.
$$
If 
\begin{equation*}
|f(p)| < \frac{2R}{\Theta (c+  n^{9/2} \|H_f(p)\|_\infty \Theta R)} =:B_2',
\end{equation*}
then the hypersurface of equation $f=0$ crosses ${\bf C}(p)$  neglecting order ${\rm O}(R^2)$ contributions.
}
\bigskip

\noindent{\bf{Remark.}} More accurate and general  estimates \cite{TBS}, when specialized to the case of hypersurfaces and $\infty$-norm,  allow us to improve the bounds above, this also assuring   that the applications discussed in \cite[Section 6]{TB} still remain  meaningful. Precisely, in \cite[Theorem 3.2]{TBS}, we can in fact show that  the bound $B_2$ goes as $\frac{1}{n}$ instead of $\frac{1}{n^{5/2}}$, weakening the assumption. 
Similarly, one can shows that the bound $B_2'$ in \cite[Theorem 4.6]{TBS} goes as $\frac{1}{n^{5/2}}$ instead of $\frac{1}{n^{9/2}}$.

\small{

}

\bigskip
\bigskip

\noindent
 M. Torrente, Istituto di Matematica Applicata e Tecnologie Informatiche - Consiglio Nazionale delle Ricerche, 
Via de Marini, 6,
I-16149 Genova, Italy. e-mail {\tt torrente@dima.unige.it}

\smallskip

\noindent M.C. Beltrametti,
Dipartimento di Matematica,
Universit\`a di Genova,
Via Dodecaneso 35,
I-16146 Genova, Italy. e-mail {\tt beltrame@dima.unige.it}


\begin{thebibliography}{10}


\bibitem{TB} M. Torrente and M.C. Beltrametti, Almost  vanishing polynomials and an application to the Hough transform, J.  Algebra Appl. {\bf 13}(8), (2014), 1450057 [39 pages].

\bibitem{TBS} M. Torrente, M.C. Beltrametti and J.R. Sendra, $r$-norm bounds  and metric properties for zero loci of real analytic functions,  preprint, 2017.


\end{thebibliography}
\end{document}


\title{Almost vanishing polynomials and an  application \\  to the Hough transform
\footnote{2010
{\em Mathematics Subject Classification}. Primary  26C10, 15A60, 14Q10;
Secondary 14H50, 68T10 \newline
\indent{{\em Keywords and phrases.}  Weighted matrix norm;  crossing area algorithm; Hough transform} }}
\author{ Maria-Laura Torrente and Mauro C. Beltrametti }

\date{}
\maketitle

\begin{abstract}
We consider the problem of deciding whether or not 
an affine hypersurface of equation $f=0$, where 
$f=f(x_1,\ldots,x_n)$ is a polynomial in $\reals[x_1,\ldots,x_n]$, 
crosses a bounded region $\mathcal T$ of the 
real affine space $\MA^n$. We perform a local study of the problem, 
and provide both necessary and sufficient numerical conditions 
to answer the question. Our conditions are based 
on  the evaluation of $f$ at a point $p \in \mathcal T$, 
and derive from the analysis of the differential geometric properties of the 
hypersurface $z=f(x_1,\ldots,x_n)$ at $p$.
We  discuss an application of our results in the context of  the Hough transform,
 a pattern recognition technique  for the automated recognition of curves in images.

\end{abstract}

\maketitle

\section*{Introduction}

We study conditions in order to decide 
whether or not a given affine hypersurface intersects a bounded region. More precisely, let $f=f({\x})$ be a polynomial in $\reals[{\x}]$, where ${\x}=(x_1,\ldots,x_n)$ denotes the variables. Intuitively, one may think that the evaluation of  $f$ at a point $p$ of the real affine space ${\mathbb A}_{\x}^n$ be sufficient to determine if either or not the hypersurface of equation $f=0$ passes through a given region containing the point $p$. That is, if $f(p)$ is sufficiently ``small"  then the hypersurfaces crosses the region, while, if $f(p)$ is ``large", the hypersurface doesn't cross the region. Indeed, this is   not the case, as the following examples show.

\begin{figure}[htb]
\centering
\begin{minipage}[c]{0.45\textwidth}
\centering
\includegraphics[width=\textwidth]{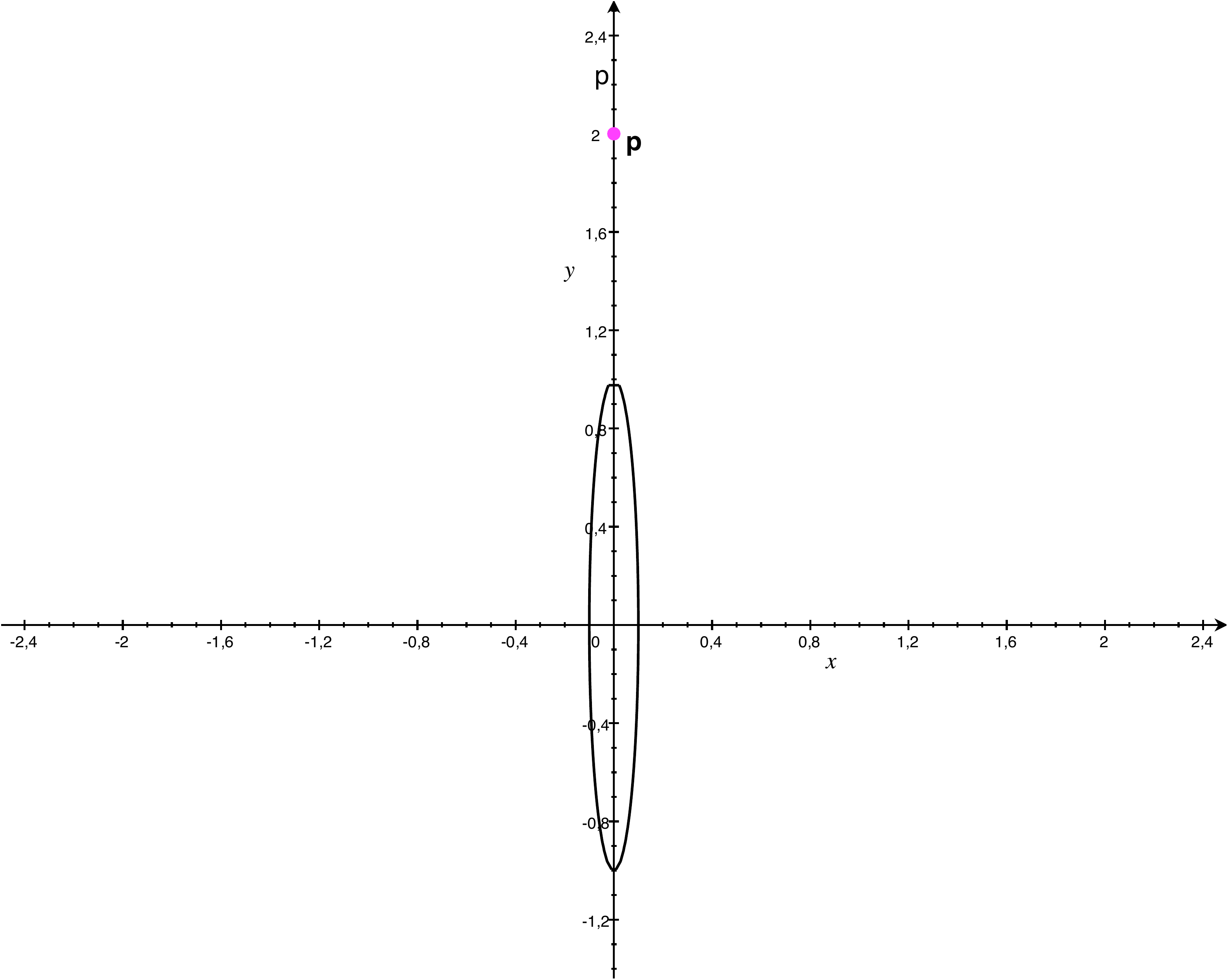}
\caption{{\small The curve $x^2+ \frac{1}{100}y^2-\frac{1}{100}=0$
and the point $p=(0,2)$.}}
\label{controEx2}
\end{minipage}%
\hspace{0.3cm}
\begin{minipage}[c]{0.48\textwidth}
\centering
\includegraphics[width=1.2\textwidth]{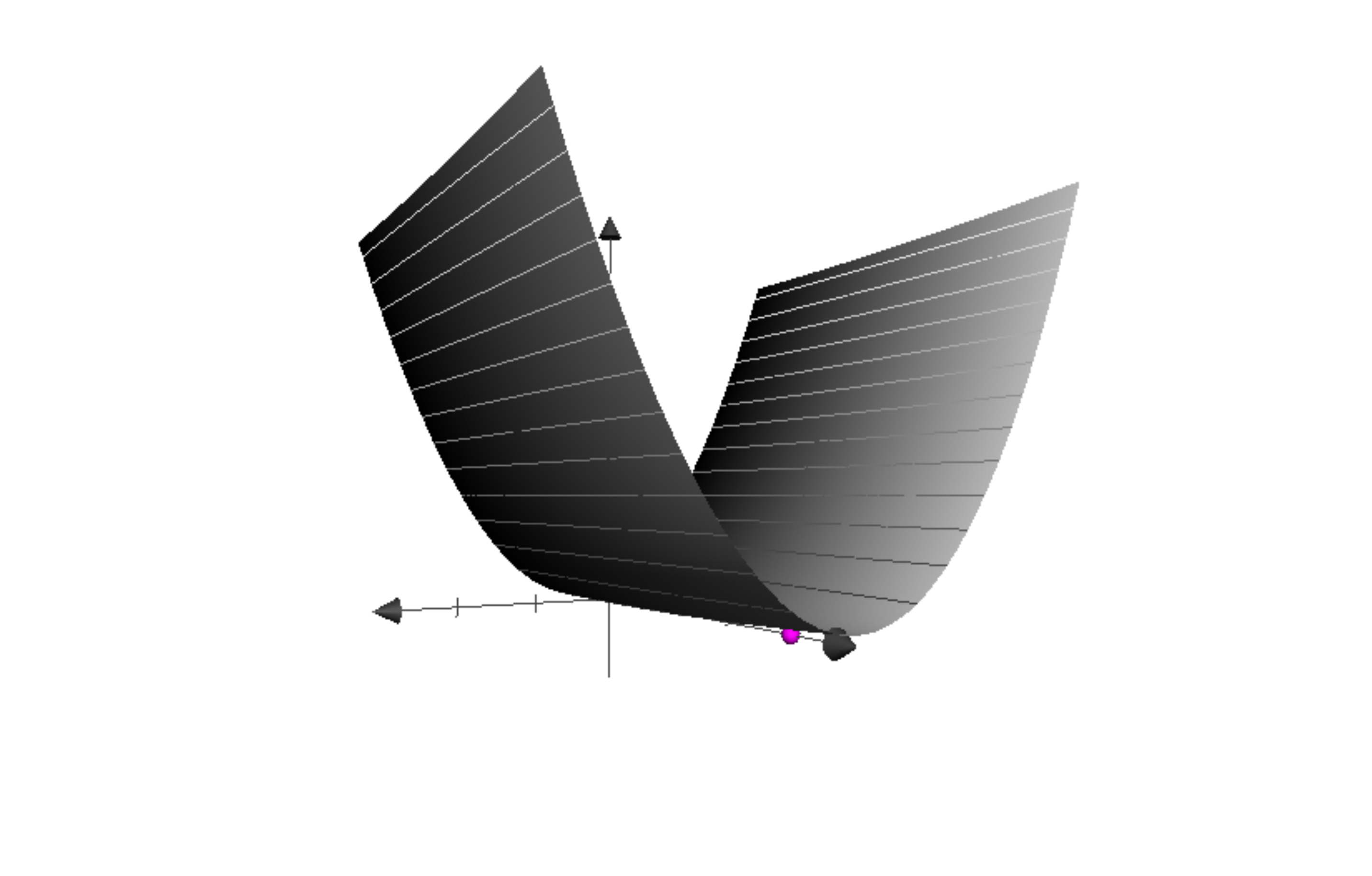}
\caption{{\small The surface $z=x^2+ \frac{1}{100}y^2-\frac{1}{100}$
and the point $p=(0,2,0)$.}}
\label{controEx2D}
\end{minipage}
\end{figure}

First, look at  the polynomial $f(x,y) =x^2+ \frac{1}{100}y^2-\frac{1}{100}
\in \R[x, y]$ and the point $p=(0,2)$. We have $|f(p)| = 0.03$ 
which  we may consider  a small evaluation. Nevertheless the point $p$ 
lies far from the curve $f=0$ (the minimal Euclidean distance of $p$ 
from points of the curve $f=0$ is~$1$,  as Figure~\ref{controEx2} shows).
The reason has to be found in the local, differential geometric properties of the  surface $z=f(x,y)$, as it is shown in Figure~\ref{controEx2D}.

\begin{figure}[htb]
\centering
\begin{minipage}[c]{0.45\textwidth}
\centering
\includegraphics[width=\textwidth]{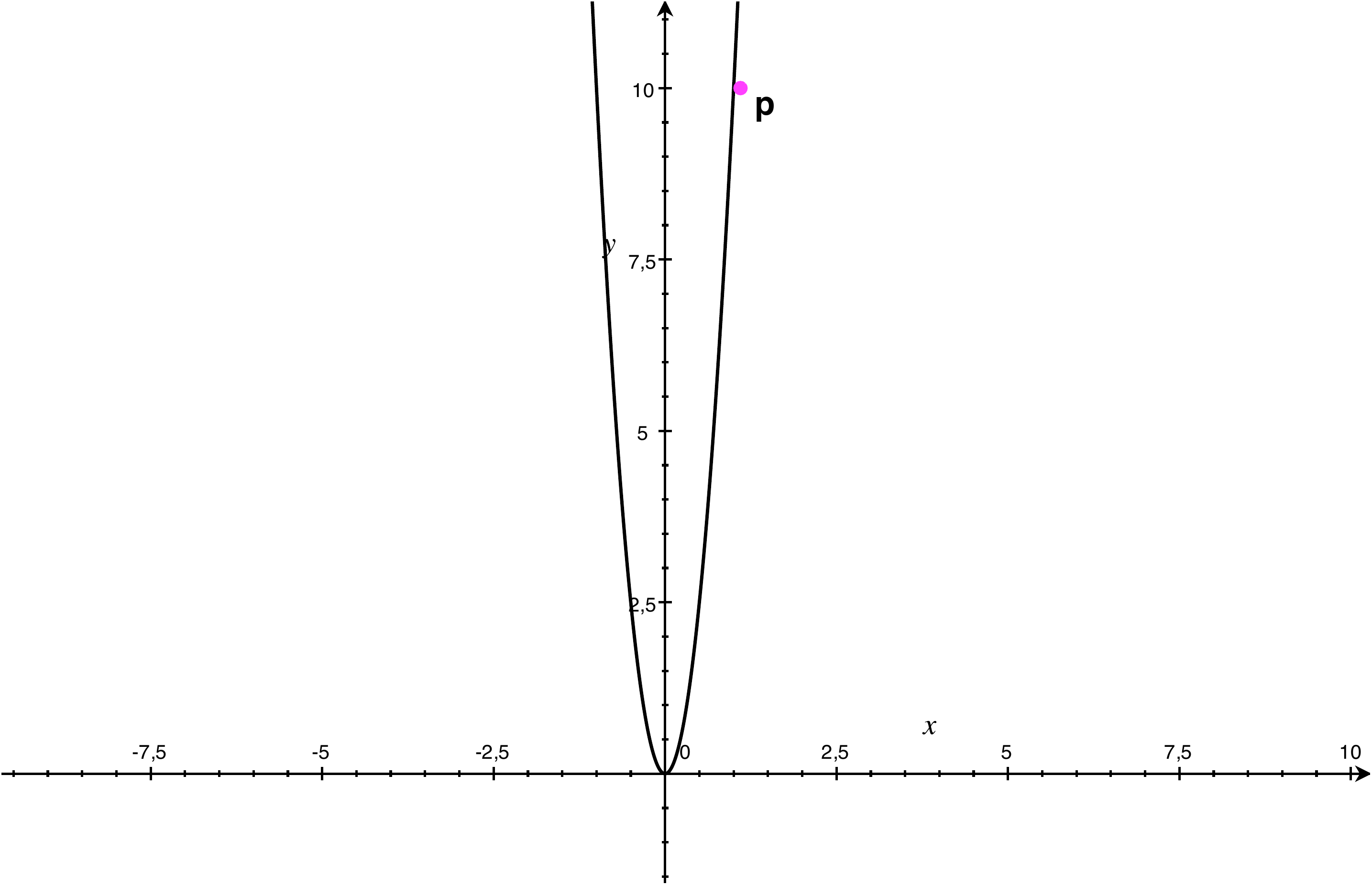}
\caption{{\small The curve $y-10x^2=0$
and the point $p=(1.1,10)$.}}
\label{controEx1}
\end{minipage}%
\hspace{0.3cm}
\begin{minipage}[c]{0.48\textwidth}
\centering
\includegraphics[width=1\textwidth]{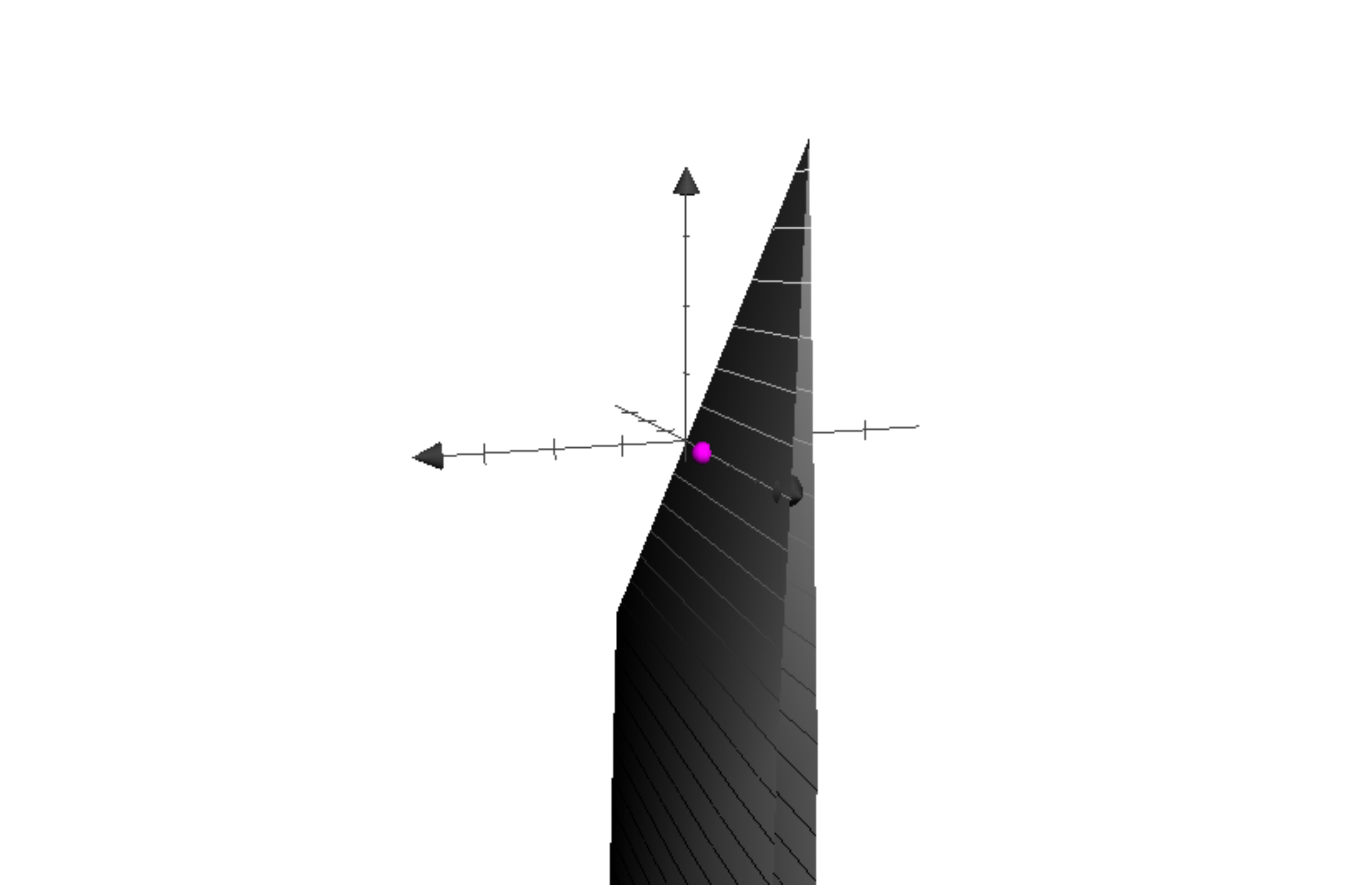}
\caption{{\small The surface $z=y-10x^2$
and the point $p=(1.1,10,0)$.}}
\label{controEx1D}
\end{minipage}
\end{figure}

Next, 
we consider the polynomial $f(x,y)=y-10x^2 \in \R[x,y]$ 
and the point $p=(1.1, 10)$. We have $|f(p)| = 2.1$ which we may think as  
a large evaluation. Nevertheless the point $p$ lies close to the curve $f=0$
(the minimal Euclidean distance of $p$ from the points of the curve $f=0$
is about~$0.1$,  as Figure~\ref{controEx1} shows).
Again, the reason has to be found in the local, differential geometric properties of the  surface $z=f(x,y)$, as it is shown in
 Figure~\ref{controEx1D}.

\smallskip
Such examples make clear a sort of ambiguity of the expression ``almost vanishing polynomials", motivating our interest to treat   the matter. The paper is organized a follows.
In Section \ref{Back} we recall the definitions of different norms and their basic properties, and a few facts about the analytic nature of polynomials.
In Section \ref{passaggioCella1} we provide necessary  numerical conditions for an affine  hypersurface $f=0$ to cross   a  bounded region containing a given  point $p$ of $\R^n$. While in Section \ref{passaggioCella2}  we provide sufficient conditions. Both the conditions (see propositions \ref{PropNonPassa} and \ref{PropPassa}) are expressed in terms of the Jacobian and the Hessian matrices, and depend on the quantities $\mathsf H$ and $\mathsf J$ (defined  in  (\ref{gamma}) and (\ref{mu}), respectively) which express exact bounds, but  not easy to be computed. 

The question to avoid the computations of $\mathsf H$ and $\mathsf J$,  by  confining ourself to find numerical conditions valid  up to a higher-order analysis, that is, up to small values of a given tolerance and disregarding higher-order contributions, becomes then natural. 
 In the case of necessary conditions we
  settle that question, up to a second-order analysis, that is, disregarding third-order contributions,
   in Proposition \ref{PropNonPassa2}.
The analogous result  for sufficient conditions  looks much harder, and it works up to a first-order analysis.   We devote the entire  Section \ref{passaggioCella2bis}   to discuss and fix this case (see Proposition \ref{PropPassa2}).

 Unfortunately, the above conditions do not fit together to give   ``if and only if " statements. However,  because of the local nature of the  results, a more accurate analysis, performed  by iteratively considering smaller subregions  of the given region 
may overcome that problem (see Remark \ref{REM}).

In Section \ref{RAL} we   summarize the results discussed above in an explanatory  algorithmic way. An implementation of the  algorithms (available at the web page {\tt http://www.dima.unige.it/$\sim$torrente/recognitionAlgorithm.cocoa5}) has been done using \cocoa 5 (see \cite{Cocoa}).

We observe that all the above mentioned conditions use the $\| \cdot\|_1$ and  $\|\cdot \|_\infty$ norms, instead of the perhaps more popular $\| \cdot \|_2$ norm. The reason for that consists in the application we have in mind,  a true guideline of our work, we deal with   
in Section \ref{RecAl}. We discuss there an explicit way of how using 
 our results in a case of special interest, 
that is,  the Hough transform technique. The Hough transform is a pattern recognition technique  for the automated recognition of curves in images
(see \cite{DH},  \cite{BMP} and  \cite{BR}). In this context, 
given a suitable family $\sF$  of curves in the image plane ${\mathbb A}_{\x}^2$,
 the Hough transform $\Gamma_p(\sF)$
 of a point $p$ in  ${\mathbb A}_{\x}^2$ with respect to $\sF$ is a hypersurface in an affine parameter space $\R^t$ (see Definition \ref{def1}), and  
 $\Gamma_p(\sF)$ plays the role of the affine hypersurface of equation $f=0$  we started from. The core of the recognition algorithm based on the Hough transforms is to count   how many hypersurfaces $\Gamma_p(\sF)$  cross a given cell of a suitable discretization of the parameter space  $\R^t$, and such a  cell is  usually defined in   $\|\cdot \|_\infty$ norm. The center of the cell which counts the maximum number of crossings is then  used to detect the recognized curve.
A relevant issue is  to validate  the outputs of our algorithm  by comparing with pattern recognition techniques exploited, e.g.,  in \cite{etal}; to this aim we provide  explicit examples which give some hopeful evidence.  A novelty here is that our approach works for any number $t$ of parameters, while, as far as we know,  the Hough transform recognition technique typically  requires  three parameters at most.

\smallskip
Special thanks are due  to A. M. Massone, 
C. Campi and  A. Perasso for  several useful discussions and  for making available experimental data we use in examples \ref{3C} and \ref{EC}. 
We are grateful to L. Robbiano for  suggesting  the possibility  to apply almost-vanishing polynomial techniques to the  Hough transform, 
and for helpful advices.  We would also like to thank J. R. Sendra for making us aware of  \cite{Sendra}, where very interesting related topics are treated.

\section{Background material}\label{Back}
\addtocounter{subsection}{1}\setcounter{theorem}{0}

In this section  we recall basic  definitions and concepts from numerical algebra and some analytic properties of polynomials systematically 
used throughout  the paper.
 We start  recalling the 
definition of different norms  on the space of matrices  and their  basic properties (for proofs and more details we refer  to \cite{Matrix}  and  \cite{TB97}).

For  $m$, $n$  positive integers, we  let ${\rm Mat}_{m \times n}(\R)$
be the set of~$m \times n$ matrices with entries in $\R$; 
if $m=n$ we simply write ${\rm Mat}_{n}(\R)$.
For any  $M \in {\rm Mat}_{m \times n}(\R)$, we will denote by $M^t$ its transpose.

\begin{definition*}\label{defnorms}
Let $v$ be  an  element of   $ {\rm Mat}_{n \times 1}(\R)$
and let $r \ge 1$ be a real number. Set $v^t:=(v_1,\ldots,v_n)$.
The {\em $r$-norm}\footnote{We will only use matrix norms; however, let us mention that in the literature one also refer to  this norm as the ``$r$-norm of the vector $(v_1,\ldots,v_n)$ in $\reals^n$".} of $v$ is defined 
by  $$\|v\|_r := \left( \sum_{i=1}^n |v_i|^r \right)^{\frac{1}{r}}.$$ In particular,
 if $r=1$, we get the expression
$
\|v\|_1 = \sum_{i=1}^n |v_i|.
$
 If $r=2$ we get the well-known {\em Euclidean norm}
$
\|v\|_2 = \big({\sum_{i=1}^n |v_i|^2}\big)^{1/2}.
$ 
While, if $r \rightarrow \infty$, the $r$-norm approaches 
the $ \infty$-{\em norm} defined by 
$
\|v\|_\infty: = \max_{i=1, \ldots, n} \big\{|v_i|\big\}.
$
\end{definition*}

\begin{prop}{\bf (Holder's inequality)} \label{propHolder}
Let  either $ q>r \ge 1$ be     real numbers such that 
$\frac{1}{q} + \frac{1}{r}=1$ or $q=\infty$ and $r=1$. 
Then for each $v, w \in  {\rm Mat}_{n \times 1}(\R)$ we have the inequality
$$
|v^t w| \le \|v\|_q \|w\|_r.
$$
\end{prop}

Apart from $r$-norms, other very useful norms 
on $ {\rm Mat}_{n \times 1}(\R)$
are  the weighted $r$-norms, in which case each component of any column vector 
$v \in  {\rm Mat}_{n \times 1}(\R)$ is rescaled according to a given weight  (see \cite{TB97}). 

\begin{definition*}\label{weightedNorm}
Let $W$ be a positive diagonal matrix in ${\rm Mat}_n(\R)$
and let $\|\cdot\|$ be any norm on $ {\rm Mat}_{n \times 1}(\R)$.
The $W$-{\em weighted norm} on  ${\rm Mat}_{n \times 1}(\R)$  is then defined  
by the formula
$$
\|v\|_W := \|Wv\|.
$$
where $v\in  {\rm Mat}_{n \times 1}(\R)$  and $W v$ denotes the usual  product of matrices. 
If either $\|\cdot\|=\|\cdot\|_r$ or $\|\cdot\|=\|\cdot\|_\infty$, we usually write
$$\|\cdot\|_{r\;  W}=:\|\cdot\|_{r,W}\;\; \;{\rm  and} \;\;\;  \|\cdot\|_{\infty \; W}=:\|\cdot\|_{\infty,W}.$$
\end{definition*}
\smallskip

The geometry of the closed unit ball (the set of all the vectors of norm 
less or equal to~$1$) clearly depends  on the norm we use. Here we consider the case of weighted $r$-norms.

\begin{definition*}\label{defCell}
Let $W$ be a positive diagonal matrix in ${\rm Mat}_n(\R)$, 
let either $r \ge 1$ be a positive real number or $r=\infty$, and 
let $p$ be a point of $\R^n$. 
The {\em $(r,W)$-unit ball centered at~$p$}, 
denoted by ${\bf B}_{r,W}(p)$, is the closed convex set defined as 
$$
{\bf B}_{r,W}(p) = \{\x \in \R^n \;\;{\rm such\; that}\;\; \|(\x-p)^t\|_{r, W} \le 1\}.
$$
For simplicity of notation, and when no confusion will arise,
the $(r,W)$-unit ball centered at $p$ will be simply called
{\em unit ball centered at $p$} and denoted by~${\bf B}(p)$.
\end{definition*}

Now, we  recall the definition of matrix norm   on $ {\rm Mat}_{m \times n}(\R)$ induced by
a given $r$-norm on $ {\rm Mat}_{n \times 1}(\R)$ and some basic facts about them.

\begin{definition*}\label{matrixNorm}
Let $M=(m_{ij})$ be a matrix in ${ \rm Mat}_{m \times n}(\R)$. 
The $r$-{\em matrix norm} is the norm on ${\rm Mat}_{m \times n}(\R) $   
induced by the $r$-norm on  ${\rm Mat}_{n \times 1}(\R)$,   and defined by the formula
$$\|M\|_r := \max_{\|v\|_r=1} \| Mv\|_r,$$  where $v\in  {\rm Mat}_{n \times 1}(\R)$.
In particular,  one has
$
\|M\|_1 = \max_{j=1,\ldots,n}\Big\{\sum_{i=1,\ldots,m} |m_{ij}| \Big\}
$ for $r=1$.
If $r=2$, denoting by $\lambda_i(\cdot)$  the $i$-th eigenvalue,   we have 
$
\|M\|_2 = \big(\max_{i=1,\ldots,n} \lambda_i(M^t M)\big)^{1/2}.
$
While,  if $r \rightarrow \infty$, the matrix $r$-norm approaches 
the $\infty$-{\em matrix  norm} defined by 
$
\|M\|_\infty: = \max_{i=1,\ldots,m}\Big\{\sum_{j=1,\ldots,n} |m_{ij}|\Big \}.
$
\end{definition*}

\begin{rem*}\label{65} Let us point out  the following consequence of Definition \ref{matrixNorm},  every so often  used throughout the paper.
 For any vector $w\in {\rm Mat}_{1\times n}(\R)$,  one has
 $$\| w\|_1= \| w^t\|_\infty.$$
\end{rem*}

Next, let us introduce one more matrix norm we need.

\begin{definition*}\label{defMaxNorm}
Let $M=(m_{ij})$ be a matrix in ${\rm Mat}_{m \times n}(\R)$. 
The {\em {\tt max}-norm} is the norm on ${\rm Mat}_{m \times n}(\R) $   
defined by 
$$
\|M\|_{\tt max} := \max_{i,j} |m_{ij}|.
$$  
\end{definition*}

Note that, via the natural identification  of $ \textrm{Mat}_n(\R)$ with $ \textrm{Mat}_{n^2 \times 1}(\R)$, 
a matrix $M\in  {\rm Mat}_n(\R)$ can be viewed as  an element, that we will denoted by $M^{({\bf v})}$ to avoid  confusion, of ${\rm Mat}_{n^2 \times 1}(\R) $. Moreover,
\begin{equation}\label{equivalent}
\| M\|_{\tt max}=\| M^{(\bf v)}\|_{ \infty}.
\end{equation}

We finally recall some useful relations between the norms introduced above, and the  sub-multiplicative property  (see \cite[\S 2.3.1]{Matrix}).

\begin{prop}\label{propMatEquivalence} 
For each $M\in {\rm Mat}_{m \times n}(\R)$ the following inequalities hold true:
\begin{enumerate}
\em\item\em \; $\frac{1}{\sqrt{n}} \|M\|_\infty \le \|M\|_2 \le \sqrt{m} \|M\|_\infty$.
\em\item\em \;  $\frac{1}{\sqrt{m}} \|M\|_1 \le \|M\|_2 \le \sqrt{n} \|M\|_1$.
\em\item\em \; $ \|M\|_{\tt max} \le \|M\|_2 \le \sqrt{mn} \|M\|_{\tt max}$.
\end{enumerate}
\end{prop}

\begin{prop}\label{subm}
Let either $r \ge 1$ be a positive real number or $r=\infty$. The $r$-matrix norm
 induced by the $r$-norm is a 
{\rm sub-multiplicative} norm, that is, for each $A \in {\rm Mat}_{m\times n}(\R)$,  $B\in {\rm Mat}_{n \times t}(\R)$, one  has 
$$
\|AB\|_r \le \|A\|_r \|B\|_r.
$$
In particular,  for each $A \in {\rm Mat}_{m\times n}(\R)$ 
and $v \in  {\rm Mat}_{n \times 1}(\R)$, one  has 
$ \|Av\|_r \le \|A\|_r \|v\|_r $.
\end{prop}

Next, we collect  some definitions 
and basic facts of analytic nature about polynomials.
The notation is borrowed from \cite{RK}.
In particular, we let $x_1,\ldots,x_n$ be indeterminates
and most of the times  we use for simplicity the notation 
$\x = (x_1,\ldots,x_n)$. 
The multivariate polynomial ring $\R[\x]=\R[x_1,\ldots,x_n]$
is denoted by~$P$.
Given $\alfa=(\alpha_1,\ldots,\alpha_n) \in \N^n$, 
we denote by~$| \alfa |$ the number $\alpha_1+\cdots+\alpha_n$, 
by $\alfa!$ the number~$\alpha_1! \ldots \alpha_n!$,  
by $\x^\alfa$ the power product $x_1^{\alpha_1} \ldots x_n^{\alpha_n}$,
and by $\frac{\partial^\alfa f}{\partial \x^\alfa} := 
\frac{\partial^{|\alfa|} f}{\partial x_1^{\alpha_1} \ldots \partial x_n^{\alpha_n}}$
the ${\bf \alfa}$-partial derivative of a polynomial $f=f(\x) \in P$. 

Moreover, following the standard notation,
we denote by $\Jac_f(\x):=\left( \frac{\partial f}{\partial x_1}, \ldots, \frac{\partial f}{\partial x_n}\right)$
 the {\em Jacobian} (or {\em gradient}) of~$f$,  and by $H_f(\x):= \left(\frac{\partial^2f}{\partial x_i\partial x_j}\right)_{i,j=1,\ldots,n}$ the $n\times n$ symmetric   {\em Hessian matrix} of~$f$.

\begin{definition*}\label{defTaylor}
Let~$p$ be a point of~$\R^n$ and 
let $f=f(\x)$ be a polynomial in~$P$. 
 Let $k$ be a non-negative integer.
The {\em $k$-th Taylor polynomial~$p_k(\x)$} and the {\em $k$-th remainder 
term}~$R_k(\x)$ of~$f(\x)$ at~$p$ are defined respectively as
$$p_k(\x) = \sum_{|\alfa| \le k} \frac{1}{\alfa!} 
\frac{\partial^\alfa f}{\partial \x^\alfa}(p) (\x-p)^\alfa \;\;\;\;{\rm and} \;\;\;\;
R_k(\x) = \sum_{|\alfa| > k} \frac{1}{\alfa!} 
\frac{\partial^\alfa f}{\partial \x^\alfa}(p) (\x-p)^\alfa,$$ so  that 
 the polynomial  $f(\x)$ can always be expressed  as $f(\x)=p_k(\x) + R_k(\x)$.
\end{definition*}

We use the following formulation of Taylor's theorem. We recall  that, given a real value $\eta\ll 1$ and  a real function $\omega:\R^n\to \R$,  we write $\omega(\x)={\rm O}(\eta^m)$, $m\in \N$, to mean that $\frac{\omega(\x)}{\eta^m}$ is bounded near the origin.

\begin{prop}{\bf (Taylor's theorem)} \label{propTaylor}
Let $k$ be a non-negative integer, let $p$ be a point of~$\R^n$ 
and let $f(\x)$ be a polynomial in~$P$. Then:
\begin{enumerate}
\em\item\em
 For each $\alfa \in \N^n$ such that $|\alfa|=k$ there exists a polynomial
$h_\alfa(\x)$  in~$P$ such that the $k$-th remainder term $R_k(\x)$ of~$f(\x)$ at~$p$ 
can be expressed as 
$$
R_k(\x) = \sum_{|\alfa|=k} h_\alfa(\x) (\x-p)^\alfa\;\;\;\;
 and\;\;\;\; \lim_{\x \rightarrow p} h_\alfa(\x) =0.$$
 In particular, $R_k(\x)={\rm O}(\|\x-p\|^{k+1})$ for any norm $\| \cdot \|$.
\em\item\em For every point $q \in \R^n$ there exists 
a point $\xi \in \R^n$ of the line segment from~$p$ to~$q$ such that the evaluation of
the $k$-th remainder term $R_k(\x)$ at $q$  is
$$
R_k(q) = \sum_{|\alfa|=k+1} \frac{1}{\alfa!}  \frac{\partial^\alfa f}{\partial \x^\alfa}(\xi)(q-p)^\alfa.
$$
\end{enumerate}
\end{prop}

The expression in Proposition \ref{propTaylor}(2) is known as the {\em Lagrange form} 
of the remainder.

Finally, we discuss a generalization of the Mean Value Theorem 
for  the case  we need of $\infty$-matrix norm.
The content of this section may be known  to experts and it is based on the following version 
of the Mean Value Theorem for vector valued real functions (see the proof presented in  \cite{HN}).
We nevertheless include  details for lack of reference.

\begin{prop}\label{propMeanValue}
Let $U \subseteq \R^n$  be a convex open set and let $p \in U$.
Let $\phi\colon U \rightarrow \R^m$ be a differentiable vector valued function on $U$
and denote by $D \phi(\x)$ the $m \times n$ matrix of first order derivatives 
of each component of $\phi$, that is,
\begin{eqnarray*}
D \phi(\x) = \left (
\begin{array}{ccc}
\frac{\partial \phi_1}{\partial x_1} & \ldots & \frac{\partial \phi_1}{\partial x_n}\\
\vdots & & \vdots\\
\frac{\partial \phi_m}{\partial x_1} & \ldots & \frac{\partial \phi_m}{\partial x_n}
\end{array}\right )
\end{eqnarray*}
Let either $r \ge 1$ be a real number or $r= \infty$. Then, for  each $\x \in U$, one has
$$
\| (\phi(\x) - \phi(p))^t \|_r < \sup_{0 < \nu< 1} \|D \phi(p+\nu(\x-p))\|_r  \|(\x -p)^t\|_r.
$$ 
\end{prop}

Let $U \subseteq \R^n$  be a convex open set,  and  let $M(\x)=(m_{ij}(\x))$ be a matrix   whose entries are the evaluations at $\x\in U$ of  differentiable vector valued  functions $m_{ij}\colon U\to \R^n$. Hence, in particular,  $M(p)\in{\rm Mat}_n(\R)$  for each given point $p\in U$.
We will use the following   special case of Proposition \ref{propMeanValue}.

\begin{lemma}\label{UppBoundH} Let $U \subseteq \R^n$  be a convex open set. Fix   a point $p$ of $U$ and 
let $M(\x)=(m_{ij}(\x))$ be a matrix  as above. For each $\x \in U$, we have
$$
\|M(\x)\|_\infty <  n^2\|M(p)\|_\infty+  {\rm O}(\|(\x-p)^t\|_\infty).
$$
\end{lemma}
\proof By combining  statements $(1)$ and $(3)$ of Proposition \ref{propMatEquivalence}
 it follows that
 \begin{equation}\label{strangeIneq}
\|M(\x)\|_\infty \le n^{3/2} \|M^{(\bf v)}(\x)\|_\infty.
\end{equation}

 Consider the vector valued function
$\phi=(M^{(\bf v)})^t\colon  U \rightarrow \R^{n^2}$ defined by $\phi(\x):= (M(\x)^{(\bf v)})^t$.
Clearly, $\phi$ is differentiable on $U$, so we can 
  apply Proposition~\ref{propMeanValue} (with $r =\infty$)  to get
$$
\|M(\x)^{(\bf v)}- M(p)^{(\bf v)}\|_\infty < 
\sup_{0<\nu<1} \|D (M(p+\nu(\x-p))^{(\bf v)})^t \|_\infty \|(\x-p)^t\|_\infty   
 ={\rm O}(\|(\x-p)^t\|_\infty).
$$
Combining the previous inequality with 
$$
\big|\: \|M(\x)^{(\bf v)}\|_\infty- \|M(p)^{(\bf v)}\|_\infty \: \big| \le \|M(\x)^{(\bf v)}- M(p)^{(\bf v)}\|_\infty
$$ (a consequence of the usual triangular inequality),
we obtain 
$$\|M(\x)^{(\bf v)}\|_\infty < \|M(p)^{(\bf v)}\|_\infty + {\rm O}(\|(\x-p)^t\|_\infty)=\|M(p)\|_{\rm max} + {\rm O}(\|(\x-p)^t\|_\infty).$$
From statements $(1)$ and $(3)$ of Proposition \ref{propMatEquivalence}   we then find 
\begin{eqnarray}\label{strange2}
\|M(\x)^{(\bf v)}\|_\infty &<& \|M(p)\|_{\rm max} + {\rm O}(\|(\x-p)^t\|_\infty) \nonumber \\
&\le& \|M(p)\|_2 + {\rm O}(\|(\x-p)^t\|_\infty) \nonumber \\
&\le& \sqrt{n}\|M(p)\|_\infty + {\rm O}(\|(\x-p)^t\|_\infty).
\end{eqnarray}
Combining (\ref{strange2}) with (\ref{strangeIneq}) we are done.
\qed

\section{Necessary crossing cell conditions }\label{passaggioCella1}
\addtocounter{subsection}{1}\setcounter{theorem}{0}
In this section we  provide necessary   numerical conditions for 
an affine hypersurface  to cross a bounded region containing 
 a given point of~$\R^n$.

 We need to fix some notation. Let $\x=(x_1,\ldots,x_n)$ be indeterminates and denote by~$P$
the multivariate polynomial ring $\R[\x]=\R[x_1,\ldots,x_n]$.
Let $f=f(\x)$ be a polynomial of~$P$ and let $p=(p_1,\ldots,p_n)$ be a point of $\R^n$. 
From now on through the paper, we make the blanket assumption that {\em the zero locus $f=0$ is of dimension $n-1$}, that is, it is a (not necessarily
irreducible) hypersurface in ${\mathbb A}_{\x}^n$, even if most of the results extends to case (not of interest  for our purposes) when the locus $f=0$
 has real components of lower dimension.

Let $\varepsilon_1,\ldots,\varepsilon_n$ be positive real numbers. 
Set  $${\epsbold:=(\varepsilon_1,\ldots,\varepsilon_n)}, \;\;
\epsbold_{{\rm min}}:= \min \{\varepsilon_1,\ldots,\varepsilon_n\},\;\;
\epsbold_{{\rm max}}:= \max \{\varepsilon_1,\ldots,\varepsilon_n\},$$
and let  ${\mathcal E} \in \textrm{Mat}_n(\R)$ be the positive diagonal matrix 
with entries $1/\varepsilon_1,\ldots,1/\varepsilon_n$. We also say that $\epsbold$ is the {\em tolerance vector}.
Throughout this section we shall use the ${\mathcal E}$-weighted $\infty$-norm 
on~$\R^n$ (see Definition~\ref{weightedNorm}) and we consider 
the corresponding closed  $(\infty,\epsbold)$-unit ball~${\bf B}(p)$ centered at $p$ 
(see Definition~\ref{defCell}). 
The results of this section apply to any set 
obtained from~${\bf B}(p)$ by removing  parts of its boundary as follows.
Consider the hyperplanes
$$L_k^{\pm}: x_k=p_k\pm \varepsilon_k, \;\; k=1,\ldots,n,$$
and let  $I^{+}, I^{-}$  be (possibly and not necessarily distinct) subsets of $\{1,\ldots,n\}$.
Then  define
\begin{equation}\label{neighP}
{\bf C}_{( I^{+}, I^{-})}(p):={\bf B}(p)\setminus\bigcup_{k\in I^{+},\;k\in I^{-}}(L_k^{+}\cup L_k^{-}).\end{equation}
 We simply write ${\bf C}(p):=
{\bf C}_{(I^{+}, I^{-})}(p)$ whenever there is no matter  what the indexes $k\in I^{+}$, $k\in I^{-}$ are considered. We  refer to  ${\bf C}(p)$ as an  $(\infty,\epsbold)$-{\em unit cell  centered at $p$}.

In the following proposition we provide a necessary  condition 
on $|f(p)|$ for  an  affine hypersurface of equation $f=0$  to cross  an $(\infty,\varepsilon)$-unit cell ${\bf C}(p)$ centered at $p$. 
Such a condition is expressed in terms of the quantity (depending on   the unit ball ${\bf B}(p)$)
\begin{equation}\label{gamma}
{\mathsf H}: = \max_{\x \in {\bf B}(p)} \|H_f(\bm x)\|_\infty.
\end{equation}

\begin{prop}\label{PropNonPassa} 
Let $f=f(\x)$ be a non-costant polynomial of~$P$,
let $p$ be a point of~$\R^n$, and let ${\bf C}(p) \subseteq {\bf B}(p)$
be an $(\infty,\epsbold)$-unit cell   centered at $p$. If
$$
| f(p)| > \|\Jac_f(p)\|_1 \epsbold_{\rm max} + \frac{\mathsf H}{2} \epsbold_{\rm max}^2 =:B_1,
$$
then the hypersurface  of equation $f=0$ does not cross ${\bf C}(p).$
\end{prop}
\proof
Suppose by contradiction that the hypersurface $f=0$ crosses 
the region ${\bf C}(p)$, that is, suppose there exists a 
point $p^* \in {\bf C}(p)$ such that $f(p^*)=0$.
From the formal Taylor expansion of $f(\x)$ at $p$ (see 
Definition~\ref{defTaylor}) it follows that
$$f(\x) =p_1(\x)+R_1(\x)= f(p) + \Jac_f(p) (\x-p)^t + R_1(\x),$$ 
where $R_1(\bm x)$ is the $1$-st remainder term of~$f(\x)$ at~$p$.
If we evaluate the former expression at~$p^*$ and apply 
Taylor's theorem (see Proposition~\ref{propTaylor}(2)) we get
\begin{eqnarray}\label{TaylorEquality}
0 = f(p^*) = f(p) + \Jac_f(p) (p^*-p)^t + \frac{1}{2}(p^*-p) H_f(\xi) (p^*-p)^t,
\end{eqnarray}
where $\xi$ is a point of the line that connects $p$ to $p^*$.
Hence,
\begin{eqnarray}\label{diseq1}
|f(p)| \le  |\Jac_f(p) (p^*-p)^t| + \frac{1}{2}|(p^*-p) H_f(\xi) (p^*-p)^t|.
\end{eqnarray}

In the following computations we systematically use H\"older's inequality  and
the sub-multiplicative property of the matrix norms (see propositions \ref{propHolder} and \ref{subm}). Recalling that 
${\mathcal E}=\textrm{diag}(\frac{1}{\varepsilon_1},\ldots,\frac{1}{\varepsilon_n})$, from  the present assumption  $p^* \in {\bf C}(p)$ 
 we obtain
\begin{eqnarray}\label{diseq2}
 |\Jac_f(p) (p^*-p)^t| & = &  |\Jac_f(p) \sE^{-1} \sE (p^*-p)^t| \nonumber \\ 
 &\le&  \|\Jac_f(p) \sE^{-1}\|_1 \|\sE (p^*-p)^t\|_\infty \nonumber \\ 
 &\le& \|\Jac_f(p) \sE^{-1}\|_1 \le  \|\Jac_f(p)\|_1 \|\sE^{-1}\|_1.
\end{eqnarray}
Analogously, we find
\begin{eqnarray}\label{diseq3}
|(p^*-p) H_f(\xi) (p^*-p)^t| & = &  |(p^*-p) H_f(\xi) \sE^{-1} \sE(p^*-p)^t| \nonumber \\ 
&\le&  \|(p^*-p) H_f(\xi) \sE^{-1}\|_1 \|\sE(p^*-p)^t\|_\infty \nonumber \\ 
&\le& \|(p^*-p) \sE \sE^{-1} H_f(\xi) \sE^{-1}\|_1 \nonumber \\
&\le& \|(p^*-p) \sE\|_1\| \sE^{-1} H_f(\xi) \sE^{-1}\|_1 \nonumber \\
&=& \|\sE(p^*-p)^t\|_\infty \| \sE^{-1} H_f(\xi) \sE^{-1}\|_1\le \|\sE^{-1}\|_1^2 \| H_f(\xi)\|_1.
\end{eqnarray}
Noting that $\|\sE^{-1}\|_1 =\epsbold_{\rm max}$ and
combining  with inequalities 
(\ref{diseq1}), (\ref{diseq2}), (\ref{diseq3}), we get
\begin{eqnarray*}
|f(p)| \le \|\Jac_f(p)\|_1 \epsbold_{\rm max} + \frac{1}{2}  \|H_f(\xi)\|_1 \epsbold_{\rm max}^2.
\end{eqnarray*}
Since $H_f(\xi)$ is a symmetric matrix, one has $\|H_f(\xi)\|_1=\|H_f(\xi)\|_\infty$. Furthermore,
 by definition, $\|H_f(\xi)\|_\infty \le {\mathsf H}$.
Thus the assertion  follows.
\qed

\begin{example*}
In $P=\R[x,y]$ we consider the polynomial
$f(x,y)=x^2+ \frac{1}{100}y^2-\frac{1}{100}$ 
and the point $p=(0,2)$,   the first  example in the Introduction.
We let $\epsbold =(0.05, 0.1)$, so
$\epsbold_{\rm max}=0.1$.
Let  ${\bf B}(p)$ be
  the $(\infty,\epsbold)$-unit cell   centered at $p$.
We have
\begin{eqnarray*}
\Jac_f(x,y)= \left(2x, \frac{1}{50}y \right)\quad \textrm{ and } \quad 
H_f(x,y)=\left ( \begin{array}{ll}
2 & 0\\
0 & \frac{1}{50}
\end{array}\right), 
\end{eqnarray*}
so that $\|\Jac_f(p)\|_1 =\frac{1}{25}$ and 
$\|H_f(x,y)\|_\infty = 2$ for each $(x,y)$. Therefore ${\mathsf H}=2$.
The magnitude $|f(p)|$ of $f$ at $p$ is $|f(p)|=0.03$, which is strictly greater than  
 the bound $B_1=0.014$ provided in Proposition~\ref{PropNonPassa}. 
We then  conclude that the curve $f=0$ does not cross ${\bf B}(p)$.
\end{example*}

\begin{example*}\label{ex1}
In $P=\R[x,y]$ we consider the polynomial $f(x,y) =4x^2+y^2-4x$.
Then
\begin{eqnarray*}
\Jac_f(x,y)= (8x-4, 2y) \quad \textrm{ and } \quad 
H_f(x,y)=\left ( \begin{array}{ll}
8 & 0\\
0 & 2
\end{array}\right). 
\end{eqnarray*}
In particular $\|H_f(x,y)\|_\infty = 8$ for each $(x,y)$, so that  ${\mathsf H}=8$.
We let $\epsbold=(0.1,0.1)$, whence $\epsbold_{\rm max}=0.1$. 
First, consider the point $p_1=(\frac{1}{4}, \frac{1}{2})$ and 
the $(\infty,\epsbold)$-unit cell  $ {\bf B}(p_1)$.
We have $|f(p_1)|=\frac{1}{2}$, which is strictly greater than  
 the bound $B_1 = 0.24$ provided in Proposition~\ref{PropNonPassa}. 
We conclude that the curve $f=0$ does not cross ${\bf B}(p_1)$.
Next,   consider the point $p_2=(\frac{1}{5}, \frac{3}{4})$ and the unit cell  ${\bf B}(p_2)$.
We have $|f(p_2)|=0.0775$, which is strictly less than
$B_1 = 0.28$. Therefore  the assumptions of Proposition~\ref{PropNonPassa}  are not satisfied. 
\end{example*}

 \begin{example*}\label{Exgamma}
In the particular case of a degree $2$ polynomial $f(\x)$ of $P$
the quantity ${\mathsf H}= \max_{\x \in {\bf B}(p)} \|H_f(\bm x)\|_\infty$ can be easily computed.

Let  $f_2(\x)=\sum_{i \le j}c_{ij} x_ix_j$  be the homogeneous 
component of~$f(\x)$ of degree $2$.
Let $p$ be a point of $\R^n$ and let
$ {\bf B}(p)$
be the $(\infty,\epsbold)$-unit ball  centered at $p$.
Then 
$$
{\mathsf H}=\max_{i=1, \ldots, n} \left (2 |c_{ii}| +
 \sum_{j \neq i} |c_{ij}| \right).
$$
To see this, recall  Definition~\ref{matrixNorm} and note that 
$H_f(\bm x)=(h_{ij})$ is the symmetric matrix with entries $h_{ij}=c_{ij}$ if $i\neq j$ and $h_{ij}=2c_{ii}$ for $i=j$.
\end{example*}

\smallskip

It is natural to ask for a statement analogous  to Proposition \ref{PropNonPassa} which  avoids the computation of   the quantity ${\mathsf H}$,  and provides  a  non-crossing cell condition simply  evaluating  the Hessian matrix at a  given point.  The following statement settles the question  for small values of (the components of) the tolerance vector $\epsbold$ and disregarding  third-order contributions. 
To this purpose, for the rest of this section, we assume $\epsbold_{\rm max}\ll 1$.

\begin{prop}\label{PropNonPassa2}
Let $f(\x)$ be a degree $\geq 2$   polynomial of~$P$.
Let $p$ be a point of $\R^n$ and let
${\bf C}(p) \subseteq {\bf B}(p)$
be an $(\infty,\epsbold)$-unit cell  centered at $p$. If
$$
|f(p)| > \|{\rm Jac}_f(p)\|_1 \epsbold_{\rm max} + \frac{1}{2}  \|H_f(p)\|_\infty \epsbold_{\rm max}^2:=B_1',
$$
then  the hypersurface of equation $f=0$ does not cross ${\bf C}(p)$ neglecting contributions of order ${\rm O}(\epsbold_{\rm max} ^3)$.
\end{prop}
\proof
Suppose by contradiction that the hypersurface $f=0$ crosses  the cell
${\bf C}(p)$, that is,  suppose there exists a 
point $p^* \in {\bf C}(p)$ such that  $f(p^*)=0$.
From the Taylor expansion of $f(\x)$ at $p$ (see 
Definition~\ref{defTaylor}) it follows that
$$f(\x) =p_2(\x)+R_2(\x)= f(p) + \Jac_f(p) (\x-p)^t + \frac{1}{2}(\x-p) H_f(p) (\x-p)^t 
+R_2(\bm x),$$ where $R_2(\bm x)$ is the 
$2$-nd remainder term of~$f(\x)$ at~$p$.
Since $p^*\in {\bf C}(p)$,  evaluating the former expression at~$p^*$ and applying 
Taylor's theorem (see Proposition~\ref{propTaylor}(1)),  we obtain 
  $R_2(p^*)={\rm O}(\|(p^*-p)^t\|_{\infty}^3)$. Thus   we find
$$
0 = f(p^*) = f(p) + \Jac_f(p) (p^*-p)^t + \frac{1}{2}(p^*-p) H_f(p) (p^*-p)^t+{\rm O}(\|(p^*-p)^t\|_{\infty}^3).
$$
Up to the term ${\rm O}(\|(p^*-p)^t\|_{\infty}^3)$, the previous expression differs from~(\ref{TaylorEquality}) only because  
$H_f(\xi)$ has been here replaced by $H_f(p)$. For this reason, and recalling that
$$\|(p^*-p)^t\|_{\infty}\leq \|\sE^{-1}\|_{\infty}\|\sE (p^*-p)^t \|_{\infty}\leq \|\sE^{-1}\|_{\infty}=\epsbold_{\rm max},$$
the same argument
 as in the proof of Proposition~\ref{PropNonPassa}
applies to give
$$ |f(p)| \leq \|{\rm Jac}_f(p)\|_1 \epsbold_{\rm max} + \frac{1}{2}  \|H_f(p)\|_\infty \epsbold_{\rm max}^2+{\rm O}(\epsbold_{\rm max}^3).$$ Thus, up to a tolerance of ${\rm O}(\epsbold_{\rm max}^3)$, we conclude that $|f(p)| \leq B_1'$. This completes the proof.
 \qed

In the following example we compare the conditions
provided in proposition~\ref{PropNonPassa} and~\ref{PropNonPassa2}. 
The bound $B_1$ is precise but harder than $B_1'$ to be computed. On the other hand, 
the bound $B_1'$, since it is effected by  a second-order 
error analysis, it is only reliable for small values of the 
tolerance vector $\epsbold$.

\begin{example*}
We consider the polynomial $f=y -10x^4 \in \mathbb R[x,y]$ 
and the point $p=(0,1)$. We let $\epsbold =(0.75,0.75)$ and 
 let $ {\bf B}(p)$
be  the $(\infty,\epsbold)$-unit ball  centered at $p$.
We observe that the plane curve of equation  $f=0$ crosses the 
neighborhood ${\bf B}(p)$. In order to verify the conditions
of propositions~\ref{PropNonPassa} and~\ref{PropNonPassa2}, evaluating    $f$ at $p$ gives $|f(p)| = 1$. Moreover,
\begin{eqnarray*}
\Jac_f(x,y)= (-40x^3, 1 ) 
\quad \textrm{ and } \quad 
H_f(x,y)=\left ( \begin{array}{cl}
-120x^2 & 0\\
0 & 0
\end{array}\right). 
\end{eqnarray*}
We have $\|\Jac_f(p)\|_1 =1$, $\|H_f(p)\|_\infty =0$,
${\mathsf H} = \max_{(x,y) \in {\bf  B}(p)}\|H_f(x,y)\|_\infty =\frac{135}{2}$. 
Consequently, the bounds $B_1$ and $B_1'$ are
$$
B_1 =\frac{1263}{64}\approx 19.7 \quad {\rm and } \quad B_1' = 0.75.
$$
The condition of Proposition~\ref{PropNonPassa} is not satisfied, 
so no conclusion can be made on the intersection of the curve of equation $f=0$
with ${\bf B} (p)$. On the other hand, the condition 
of Proposition~\ref{PropNonPassa2} is verified, so, up to a second-order 
error analysis, we would  wrongly conclude that the curve $f=0$ does not cross~${\bf  B}(p)$ (in this case, the value $0.75$ of the components of $\epsbold $ is not small enough).
\end{example*}

We end this section  by comparing  the bounds  provided by propositions \ref{PropNonPassa} and \ref{PropNonPassa2}.

\begin{lemma}\label{PropB1B1'} Let $B_1$, $B_1'$ be the bounds as above. Thus
$B_1 \ge B_1'$.
\end{lemma}
\proof
It follows from the definition of $B_1$ and $B_1'$,
and from the obvious inequality ${\mathsf H}  = \max_{\x \in {\bf  B}(p)} \|H_f(\x)\|_\infty
\ge \|H_f(p)\|_\infty$.
\qed

\section{Sufficient crossing cell conditions I}\label{passaggioCella2}
\addtocounter{subsection}{1}\setcounter{theorem}{0}

In this section we provide sufficient  numerical conditions 
for  an  affine hypersurface to cross  a bounded region containing 
a given point  of~$\R^n$. 

 For each $\x=(x_1,\ldots,x_n)$ such that ${\rm Jac}_f(\x)$ is not zero, we consider the  {\em pseudo-inverse matrix} of ${\rm Jac}_f(\x)$, defined by
$$
{\rm Jac}_f^\dagger(\x):={\rm Jac}_f(\x)^t\big({\rm Jac}_f(\x) {\rm Jac}_f(\x)^t\big)^{-1} \Big(=\frac{{\rm Jac}_f(\x)^t}{\|{\rm Jac}_f(\x)\|_2^2}\Big).$$ 
Note  that ${\rm Jac}_f^\dagger(\x)$ is the right inverse of $ {\rm Jac}_f(\x)$, that is, 
${\rm Jac}_f(\x){\rm Jac}_f^\dagger(\x)=1$.

For any    positive real number $R$,  set
$$
\mathcal D(p, R):=\{\x \in \R^n \;\;{\rm such\; that}\; \; \|(\x -p)^t \|_\infty < R\}.
$$ 
Note that  $\mathcal D(p, R)\subseteq \mathcal {\bf C}(p)$ as soon as $R< \epsbold_{\rm min}$. 
Then set
\begin{equation}\label{mu}
\mathsf J: = \sup_{\x \in \mathcal D(p, R)} \|{\rm Jac}_f^\dagger(\x)\|_\infty.
\end{equation}

We need the following technical lemma.

\begin{lemma}\label{Full1} 
Let $f=f(\x)$ be a degree $\geq 2$ polynomial of~$P$ and 
let $p$ be a point of $\R^n$ such  that the Jacobian  $\Jac_f(p)$ is nontrivial.
Let $R$ be a positive real number such that $R < \epsbold_{{\rm min}}$.
If $ R<\frac{\|\Jac_f(p)\|_1}{\mathsf H}$, 
 then ${\rm Jac}_f(\x)$ is   nontrivial for  $\x\in  \sD(p,R)$. 
\end{lemma}
\proof
First we note that since $R < \epsbold_{\rm min}$ we have $\sD:=\sD(p,R) \subseteq {\bf B}(p)$.
Moreover, for each point $\x \in \sD$, one has  the triangular inequality  relation
\begin{equation}\label{miserve}
\big| \|\Jac_f(\x)\|_1 - \|\Jac_f(p)\|_1 \big| \le \|\Jac_f(\x) - \Jac_f(p)\|_1.
\end{equation}
 
 We consider the vector valued function $\phi = \Jac_f\colon \sD \rightarrow \R^n$;
we observe that $\phi$ is differentiable on the open convex set $\sD$
and that $D \phi(\x) = H_f(\x)$. Applying Proposition~\ref{propMeanValue} 
to $\phi$ with $r =\infty$ and using the fact that $\x \in {\bf B}(p)$, we then have
\begin{eqnarray}\label{depressione}
\|\Jac_f(\x)^t - \Jac_f(p)^t\|_\infty &<& \sup_{0<\nu<1} \|H_f(p+\nu(\x-p))\|_\infty \|(\x-p)^t\|_\infty 
\nonumber\\
&<& \sup_{0<\nu<1} \|H_f(p+\nu(\x-p))\|_\infty R.
\end{eqnarray}
Recalling  definition (\ref{gamma}) of the quantity $\mathsf H$, we get
$
\|\Jac_f(\x)^t - \Jac_f(p)^t\|_\infty < \mathsf H R.
$
Since $\|\Jac_f(\x)^t - \Jac_f(p)^t\|_\infty = \|\Jac_f(\x) - \Jac_f(p)\|_1$ (see Remark \ref{65}),
from inequality (\ref{miserve}) we then find  
$\|\Jac_f(\x)\|_1 > \|\Jac_f(p)\|_1 - {\mathsf H} R >0$, 
so that  the Jacobian  $\Jac_f(\bm x)$ is nontrivial  for $\x\in \sD$.
\qed

We have the following sufficient condition to conclude that  the hypersurface of equation $f=0$  
crosses  a given unit cell ${\bf C}(p)$. 

\begin{prop}\label{PropPassa} 
Let $f=f(\x)$ be a degree $\geq 2$ polynomial of~$P$,
let $p$ be a point of $\R^n$ such that $\Jac_f(p)$ is not the zero vector, 
and let ${\bf C}(p) \subseteq {\bf B}(p)$
be an $(\infty,\epsbold)$-unit cell  centered at $p$.
Let $R$ be a positive real number such that 
$R < \min \big\{\epsbold_{\rm min}, \frac{\|\Jac_f(p)\|_1}{{\mathsf H}}\big\}$.
Set $c:=\max\{2, \sqrt{n}\}$. If 
$$
|f(p)| < \frac{2R}{\mathsf J(c+  \sqrt{n}{\mathsf H} {\mathsf J} R)}  =: B_2,
$$
then the hypersurface of equation $f=0$ crosses ${\bf C}(p)$.
\end{prop}
\proof  If $f(p)=0$ there is nothing to prove.
From Lemma \ref{Full1} we know that the Jacobian   $\Jac_f(\bm x)$ is nonzero for   $x\in \mathcal D:=\mathcal D(p, R)$.
Moreover, since $R< \epsbold_{\rm min}$, one has $\sD\subseteq \mathcal {\bf C}(p)$.

We now construct a sequence of points $\{p_k\}_{k \in \N}$ as follows. 
We let $p_0=p$ and, for each $k \ge 0$, we define 
\begin{equation}\label{DEF1}
s_k:=-\Jac_f^\dagger(p_k)  f(p_k)= - \frac{f(p_k)}{\|{\rm Jac}_f(p_k)\|_2^2} {\rm Jac}_f(p_k)^t \qquad {\rm and}\qquad
p_{k+1}:=p_k+s_k^t.
\end{equation}
The construction of the points $p_k$ draws back to the Normal Flow algorithm
(see \cite{WW}), an iterative method mainly used in homotopy  and continuation 
problems. 
Obviously, $p=p_0\in \sD$.
We prove by induction that the points  $p_k$'s all lie in $ \mathcal D$ and satisfy the
 inequality 
 \begin{equation}\label{civuole}
 |f(p_k)| <  |f(p_{k-1})| \;\;{\rm for \; each}\;\;k\geq 1.\end{equation}

\noindent{\bf Step I} (The $k=1$ case). 
From the definitions of $s_0$ and $\mathsf J$ (see (\ref{mu})) we have
$$
\|s_0\|_\infty = \|{\rm Jac}_f^\dagger(p)\|_\infty |f(p)| \le {\mathsf J} |f(p)|.
$$
Moreover, by assumption, it follows that
$|f(p)| <B_2 < \frac{2 R}{c \mathsf J} \le \frac{R}{\mathsf J}$.  Thus $\|s_0\|_\infty < R$ (since $\sD$ is defined in $\| \cdot  \|_\infty$ norm), showing that  $p_1 \in \mathcal D$.

From the formal Taylor expansion of $f(\x)$ at $p$ (see 
Definition~\ref{defTaylor}) it follows that
$$
f(\x) = f(p) + \Jac_f(p) (\x-p)^t + R_1(\x),
$$
where $R_1(\bm x)$ is the $1$-st remainder term of~$f(\x)$ at~$p$.
If we evaluate the former expression at~$p_1$ and apply 
Taylor's theorem (see Proposition~\ref{propTaylor}(2)) we get
\begin{eqnarray*}
f(p_1) = f(p) + \Jac_f(p) (p_1-p)^t + 
\frac{1}{2}(p_1-p)H_f(\xi) (p_1-p)^t,
\end{eqnarray*}
 where $\xi$ is a point of the line that connects $p$ to $p_1$.
Therefore,  by definitions (\ref{DEF1}), we get
\begin{eqnarray}\label{caldo} \nonumber
f(p_1)  &=&  f(p) + \Jac_f(p) s_0 +  \frac{1}{2} s_0^t H_f(\xi) s_0\\ \nonumber
&=&  f(p) - \frac{f(p)}{\|{\rm Jac}_f(p)\|_2^2} \Jac_f(p)  \Jac_f(p)^t +  
\frac{1}{2} \frac{|f(p)|^2}{\|{\rm Jac}_f(p)\|_2^4} \Jac_f(p) H_f(\xi) \Jac_f(p)^t\\ \nonumber
&=& f(p) - f(p) +  
\frac{1}{2} \frac{|f(p)|^2}{\|{\rm Jac}_f(p)\|_2^4} \Jac_f(p) H_f(\xi) \Jac_f(p)^t\\ 
&=&  |f(p)| \left ( 
\frac{1}{2} \frac{|f(p)|}{\|{\rm Jac}_f(p)\|_2^4} \Jac_f(p) H_f(\xi) \Jac_f(p)^t \right ).
\end{eqnarray}

Let us upper bound the absolute value of the quantity
$$Q:= \frac{1}{2} \frac{|f(p)|}{\|{\rm Jac}_f(p)\|_2^4} \Jac_f(p) H_f(\xi) \Jac_f(p)^t. $$
To this end, use   H\"older's inequality (see Proposition~\ref{propHolder}), Proposition  \ref{propMatEquivalence}(1) 
and Remark \ref{65},
 and recall   the definitions of ${\mathsf H}$, $\mathsf J$ (see (\ref{gamma}),  (\ref{mu}) respectively) to get:
\begin{eqnarray*}
|Q| &\le & \frac{1}{2} \frac{|f(p)|}{\|{\rm Jac}_f(p)\|_2^4}  \| \Jac_f(p)\|_2^2 \: \|H_f(\xi)\|_2
\le  \frac{1}{2} \frac{|f(p)|}{\|{\rm Jac}_f(p)\|_2^2}  \sqrt{n} \|H_f(\xi)\|_\infty\\
&\le & \frac{1}{2} \frac{ |f(p)|}{\|{\rm Jac}_f(p)\|_2^2} \sqrt{n} {\mathsf H} =
 \frac{1}{2} \sqrt{n}|f(p)| \frac{\|\Jac_f(p)^t\|_\infty}{\|{\rm Jac}_f(p)\|_2^2}  
 \frac{\mathsf H}{\|\Jac_f(p)^t\|_\infty} \\
 &=& \frac{1}{2} \sqrt{n}|f(p)| \|\Jac_f^\dagger(p)\|_\infty  \frac{{\mathsf H}}{\|\Jac_f(p)\|_1}   
\le  \frac{1}{2} \sqrt{n}|f(p)| \mathsf J  \frac{{\mathsf H}}{\|\Jac_f(p)\|_1}.    
\end{eqnarray*}
By the assumption on $R$   we thus obtain
$|Q|< \frac{1}{2} \sqrt{n}|f(p)| \frac{\mathsf J}{R}$.
On the other hand,  $|f(p)| < B_2 < \frac{2R}{c \mathsf J} \le  \frac{2R}{\sqrt{n} \mathsf J}$.
 Therefore  $|Q|<1$,  so that equality (\ref{caldo}) reads $|f(p_1)| <  |f(p)|$, 
   showing condition  (\ref{civuole}) for $k=1$.
   \smallskip
 
\noindent{\bf Step II} (The inductive step).  
 Suppose that the points $p,p_1, \ldots p_k$ of the sequence lie in $\mathcal D$
and that $0< |f(p_k)| < |f(p_{k-1})| <  \cdots < |f(p)|$. Hence, in particular, the points $p,p_1, \ldots p_k$ are all distinct, so that, by definition, $\|s_{i-1}\|_\infty\neq 0$ for $i=1,\ldots,k$.

First we show that $p_{k+1} \in \mathcal D$. For $i=1,\ldots,k$,  the formal Taylor expansion of $f(\x)$ at $p_{i-1}$ (see 
Definition~\ref{defTaylor})  yields
$$
f(\x) = f(p_{i-1}) + \Jac_f(p_{i-1}) (\x-p_{i-1})^t + R_1(\x),
$$
where $R_1(\bm x)$ is the $1$-st remainder term of~$f(\x)$ at~$p_{i-1}$.
If we evaluate the former expression at~$p_i$ and apply 
Taylor's theorem (see Proposition~\ref{propTaylor}(2)) we get
\begin{equation}\label{EqT1}
f(p_i) = f(p_{i-1}) + \Jac_f(p_{i-1}) (p_i-p_{i-1})^t + 
\frac{1}{2}(p_i-p_{i-1})H_f(\xi_i) (p_i-p_{i-1})^t,
\end{equation}
where $\xi_i$ is a point of the line that 
connects $p_{i-1}$ to $p_i$.
On the other hand, by definition of $s_{i-1}$ and recalling that ${\rm Jac}_f(\x){\rm Jac}_f^\dagger(\x)=1$,  we have
$\Jac_f(p_{i-1})s_{i-1}=-f(p_{i-1})$, whence 
\begin{equation}\label{EqT2}
f(p_{i-1})=    -\Jac_f(p_{i-1})s_{i-1}=   - \Jac_f(p_{i-1})(p_i-p_{i-1})^t.
\end{equation}
By combining (\ref{EqT1}) and (\ref{EqT2})  with
H\"older's inequality, we get
\begin{eqnarray}\label{idem} \nonumber 
|f(p_i)| &=&  \frac{1}{2}|(p_i-p_{i-1}) H_f(\xi_i)(p_i - p_{i-1})^t| \\ \nonumber
&\le&\frac{1}{2} \|(p_i-p_{i-1}) H_f(\xi_i)\|_1 \|(p_i - p_{i-1})^t\|_\infty \\ \nonumber
&\le& \frac{1}{2} \|p_i-p_{i-1}\|_1  \|H_f(\xi_i)\|_1 \|(p_i - p_{i-1})^t\|_\infty \\ 
&=& \frac{1}{2}  \|H_f(\xi_i)\|_1 \|(p_i - p_{i-1})^t\|_\infty^2.
\end{eqnarray}
Since the Hessian matrix is symmetric we have
$\|H_f(\xi_i)\|_1=\|H_f(\xi_i)\|_\infty$. Thus,  by definition of 
${\mathsf H}$,  the previous relation yields
\begin{eqnarray}\label{ineq1}
|f(p_i)| \le \frac{1}{2}  \|H_f(\xi_i)\|_\infty \|(p_i- p_{i-1})^t\|^2_\infty  
\leq \frac{1}{2} \mathsf H \|s_{i-1}\|_\infty^2.
\end{eqnarray}

Now,  
define $\tau_i : =\frac{\|s_i\|_\infty}{\|s_{i-1}\|_\infty}$.
Therefore inequality (\ref{ineq1}) gives 
$$
\|s_i\|_\infty = \|\Jac_f^\dagger(p_i)\|_\infty |f(p_i)|
\le {\mathsf J} |f(p_i)| \le \frac{1}{2} \mathsf J \mathsf H \|s_{i-1}\|_\infty^2.
$$
Thus
\begin{eqnarray}\label{ineqTau0}
\tau_i = \frac{\|s_i\|_\infty}{\|s_{i-1}\|_\infty} \le \frac{1}{2} \mathsf J \mathsf H \|s_{i-1}\|_\infty
\le \frac{1}{2} \mathsf J^2 \mathsf H |f(p_{i-1})|  < \frac{1}{2} \mathsf J^2 {\mathsf H} |f(p)|. 
\end{eqnarray}
Since $|f(p)| < B_2 < \frac{2}{\sqrt{n} \mathsf J^2 {\mathsf H}} \le \frac{2}{\mathsf J^2 {\mathsf H}}$,
 it must be $\tau_i<1$ by the above inequality. 
Let $\tau:=\max_{i=1, \ldots, k} \{\tau_i\}$. We bound $\|(p_{k+1}-p)^t\|_\infty$ as follows:
\begin{eqnarray*}
 \|(p_{k+1} -p)^t\|_\infty &\le& \|s_0\|_\infty + \|s_1\|_\infty + \cdots +\|s_k\|_\infty \\
 &=& \|s_0\|_\infty + \tau_1 \|s_0\|_\infty + \tau_1\tau_2 \|s_0\|_\infty +  \cdots +\tau_1 \tau_2 \ldots \tau_k \|s_0\|_\infty\\
&=& \|s_0\|_\infty (1 + \tau_1+\tau_1\tau_2 +\cdots + \tau_1\tau_2 \dots \tau_k)\\
&\leq& \|s_0\|_\infty \sum_{i=0}^k \tau^i <\|s_0\|_\infty  \sum_{i=0}^\infty \tau^i =\frac{\|s_0\|_\infty}{1- \tau}\leq   \frac{\mathsf J |f(p)|}{1- \tau}.
\end{eqnarray*}
Then, by  inequality (\ref{ineqTau0}) and  the assumption $|f(p)| < B_2$, we find
$$
\|(p_{k+1} -p)^t\|_\infty < \frac{\mathsf J |f(p)|}{1- \frac{1}{2} \mathsf J^2 {\mathsf H} |f(p)|} =
\frac{2 \mathsf J |f(p)|}{2- \mathsf J^2 {\mathsf H} |f(p)|} < R,
$$
therefore $p_{k+1} \in \mathcal D$.

Now, let us prove that $|f(p_{k+1})| < |f(p_k)|$. To this purpose we observe
that relation (\ref{caldo}) can be easily adapted to the pair of points $p_k$, $p_{k+1}$
in the form
\begin{eqnarray}\label{fpk}
f(p_{k+1}) = |f(p_k)| \left ( 
\frac{1}{2} \frac{|f(p_k)|}{\|{\rm Jac}_f(p_k)\|_2^4} \Jac_f(p_k) H_f(\xi_k) \Jac_f(p_k)^t \right ),
\end{eqnarray}
where $\xi_k$ is a point of the line connecting $p_k$ to $p_{k+1}$.
Let us upper bound the absolute value of the quantity
$$
Q_k:= \frac{1}{2} \frac{|f(p_k)|}{\|{\rm Jac}_f(p_k)\|_2^4} \Jac_f(p_k) H_f(\xi_k) \Jac_f(p_k)^t.
$$
As  previously done to upper bound the quantity $|Q|$, by using now H\"older's inequality, the first two statements  of 
Proposition \ref{propMatEquivalence}, and the definition of $\mathsf H$, we get
\begin{eqnarray}\label{ineqQk}
|Q_k| 
&\le&  \frac{1}{2} \sqrt{n} \frac{|f(p_k)|}{\|{\rm Jac}_f(p_k)\|_2^2}   \|H_f(\xi_k)\|_\infty
\le  \frac{1}{2}\sqrt{n}  \frac{|f(p_k)|}{\|{\rm Jac}_f(p_k)\|_2^2}  \mathsf H \nonumber\\
&\le&  \frac{1}{2}\sqrt{n}  \frac{|f(p_k)|}{\|{\rm Jac}_f(p_k)\|_1^2}  \mathsf H 
<  \frac{1}{2}\sqrt{n}  \frac{|f(p)|}{\|{\rm Jac}_f(p_k)\|_1^2}  \mathsf H, 
\end{eqnarray}
where the last inequality comes from the inductive hypothesis $|f(p_k)| < |f(p)|$.
Since $\Jac_f(p_k) \Jac_f(p_k)^\dagger = 1$, 
H\"older's inequality and the definition of $\mathsf J$ give
$$
1= |\Jac_f(p_k) \Jac_f(p_k)^\dagger| \le \|\Jac_f(p_k)\|_1  \|\Jac_f(p_k)^\dagger\|_\infty
\le \|\Jac_f(p_k)\|_1 \mathsf J. 
$$
Therefore inequality (\ref{ineqQk}) becomes
$|Q_k| <  \frac{1}{2}\sqrt{n}  |f(p)|  {\mathsf H} {\mathsf J}^2$.
On the other hand, $|f(p)| < B_2 < \frac{2}{\sqrt{n} \mathsf{H} \mathsf{J}^2}$.
Thus we  find $|Q_k|<1$,  so that equality (\ref{fpk})  
yields $|f(p_{k+1})| < |f(p_k)|$, as we want.
\smallskip

\noindent{\bf Step III} (Conclusion). If there exists $k\in \N$ such that $f(p_k)=0$ we are done.
 Otherwise, 
 we know from Step II  that $\tau_k:=\frac{ \|s_k\|_\infty}{ \|s_{k-1}\|_\infty}<1$ for $k\in \N$.  Then, by D'Alembert criterion, the series $\sum_{k=1}^\infty\|s_k\|_\infty$ converges, so that $\lim_{k\to \infty}\big(\sum_{i=k+1}^\infty\|s_i\|_\infty\big)=0$. Define
${p^*}^t:=p^t+\sum_{k=1}^\infty s_k$. Then, since $p_k^t=p^t+\sum_{i=1}^ks_i$, one has
$$\lim_{k\to\infty}\|(p_k-p^*)^t\|_\infty=\lim_{k\to\infty}\Big(\|\sum_{i=k+1}^\infty s_i \|_\infty\Big)\leq \lim_{k\to \infty}\Big(\sum_{i=k+1}^\infty\|s_i\|_\infty\Big)=0.$$
Thus  the sequence of points $\{p_k\}_{k \in \N}$ converges to the point $p^*$. Since the  $p_k$'s belong to $\sD$, the point $p^*$ belongs to the closure $\overline{\sD}\subseteq {\bf C}(p)$.
We also know that $\| s_k\|_\infty=\tau_1\tau_2\ldots\tau_k\|s_0\|_\infty<\tau^k\|s_0\|_\infty$, where $\tau=\sup_{k\in \N}\{\tau_k\}$. Therefore $\lim_{k\to \infty}\|s_k\|_\infty< \lim_{k\to \infty}\tau^k\|s_0\|_\infty=0$.
From  inequality~(\ref{ineq1}), 
 we then conclude that 
$|f(p^*)|=\lim_{i\to \infty} |f(p_i)|\le  \frac{1}{2} \mathsf H\lim_{i\to \infty} \|s_{i-1}\|_\infty^2=0$. This completes the proof. 
\qed

\begin{rem*}\label{lineare} (The linear case) Notation as above. It is just the case noting that, if  the polynomial $f=f(\x)$ is  linear, the bound $B_1$ of Proposition \ref{PropNonPassa} simply becomes $B_1= \|\Jac_f(p)\|_1 \epsbold_{\rm max}$. 

Concerning the results of the present section, in the  linear case, Lemma 
 \ref{Full1} holds true under the only assumption that $R<\epsbold_{{\rm min}}$.
Similarly,  Proposition \ref{PropPassa} holds true  under the only assumption that $R< \epsbold_{{\rm min}}$  as well. Moreover,  the bound $B_2$ simply becomes $B_2=\frac{R}{\mathsf J}$.  Indeed, for a linear polynomial $f$,  the same argument as in Step I of the proposition shows that the hyperplane $f=0$ crosses the cell ${\bf C}(p)$ as soon as $|f(p)|<\frac{R}{\mathsf J}$ (since equation (\ref{caldo}) yields $f(p_1)=0$ in that case).
 \end{rem*}

\begin{example*}
We consider the polynomial $f=y -10x^2$ 
and the point $p=(1.1,10)$, the second  example in the Introduction.  
We let $\epsbold =(0.13,0.13)$ and consider the
$(\infty,\epsbold)$-unit ball ${\bf B}(p)$ centered at $p$.
We have
\begin{eqnarray*}
{\rm Jac}_f(\bm x)= (-20x, 1 )\quad \textrm{ and } \quad 
{\rm Jac}_f^\dagger(\bm x)= \frac{1}{1+400x^2}
\left ( \begin{array}{c}
-20x\\ 1  \end{array}\right ),
\end{eqnarray*}
whence  $\|{\rm Jac}_f(p)\|_1 =22 >0$. Furthermore $\|H_f(\bm x)\|_\infty = 20$, 
for each $\bm x$, so that ${\mathsf H}=20$.
We let $R =0.12$ and compute $\mathsf J \approx 0.05$.  
We have that $|f(p)|=2.1$ is strictly smaller than  
$\frac{2}{\mathsf J(c+ \sqrt{n}{\mathsf H} {\mathsf J} R)} R \approx 2.21$,
which is the bound provided in Proposition~\ref{PropNonPassa}. 
We conclude that the curve of equation $f=0$ goes through the ball ${\bf B}(p)$.
\end{example*}

\begin{example*}
We consider the polynomial $f =4x^2+y^2-4x$ as in Example~\ref{ex1} and the point 
$p_2=(\frac{1}{5}, \frac{3}{4})$.
We let $\epsbold =(0.1,0.1)$ and consider the
$(\infty,\epsbold)$-unit ball   ${\bf B}(p)$ centered at $p$.
We have ${\mathsf H}=8$  and 
$\|\Jac_f(p_2)\|_1=\|(-2.4, 1.5)\|_1 = 2.4 >0$. We choose $R=0.075 < \min \{0.1, 0.3\}$.
Direct computations  show that 
$\mathsf J =\sup_{(x,y) \in \mathcal D} \|\Jac_f^\dagger(x,y)\|_1 = \frac{16}{45}$.
Therefore for the bound $B_2$ of Proposition \ref{PropPassa}  we find  
 $ B_2 = \frac{2R}{\mathsf J(c+ \sqrt{n}{\mathsf H} {\mathsf J} R)} \approx 0.18$.
Since $|f(p)|=0.0775$ is strictly smaller than $B_2$, 
by using Proposition~\ref{PropNonPassa} we conclude that the 
curve of equation $f=0$ crosses the  ball ${\bf B}(p)$. 
\end{example*}

\section{Sufficient crossing cell conditions II}\label{passaggioCella2bis}
\addtocounter{subsection}{1}\setcounter{theorem}{0}
In this section we provide sufficient  numerical conditions, working up to a first-order error analysis,  for 
 an affine hypersurface  to cross
 a bounded region containing a given point  of~$\R^n$.  

To begin with, let us  upper bound the quantity $\mathsf J$ defined in  (\ref{mu}).

\begin{prop}\label{1Luglio}
Let $f(\x)$ be a  degree $\geq 2$   polynomial of~$P$ and 
let $p$ be a point of $\R^n$. Let $R$ be a positive real number and 
suppose that the Jacobian $\Jac_f(\x)$ is nonzero for each 
$\x \in \mathcal D(p,R)$. Then 
$$
\mathsf J <  
 \|\Jac_f^\dagger(p)\|_\infty+  n^2(1+2 \sqrt{n}) \frac{ \|H_f(p)\|_\infty}{\|\Jac_f(p)\|_1^2} R +{\rm O}(R^2).
$$
\end{prop}
\proof Let $\x \in \mathcal D:= \mathcal D(p,R)$. Consider the vector valued function {$(\Jac_f^\dagger)^t \colon \sD \rightarrow \R^n$.
Since by hypothesis $\Jac_f(\x)$ has full row rank in $\mathcal D$,
it follows that  $(\Jac_f^\dagger)^t$ is differentiable on the open convex set $\mathcal D$.
We apply Proposition~\ref{propMeanValue} with $r =\infty$ to get 
\begin{equation}\label{2410}
\|\Jac_f^\dagger(\x) - \Jac_f^\dagger(p)\|_\infty < 
\sup_{0<\nu<1} \|D \Jac_f ^\dagger(p+\nu(\x-p))\|_\infty \|(\x-p)^t\|_\infty.
\end{equation}
Combining (\ref{2410})  with
$\big| \|\Jac_f^\dagger(\x)\|_\infty - \|\Jac_f^\dagger(p)\|_\infty \big| \le \|\Jac_f^\dagger(\x) - \Jac_f^\dagger(p)\|_\infty$
(the usual consequence of the triangular inequality),
 we have
\begin{eqnarray}\label{princDisugJac}
\|\Jac_f^\dagger(\x)\|_\infty < \|\Jac_f^\dagger(p)\|_\infty + 
\sup_{0<\nu<1} \|D \Jac_f^\dagger (p+\nu(\x-p))\|_\infty \|(\x-p)^t\|_\infty.
\end{eqnarray}
Applying Lemma \ref{UppBoundH} to the matrix $M(\x)= D \Jac_f^\dagger (\x)$, one has
\begin{eqnarray}\label{supDJac}
\sup_{0<\nu<1} \|D \Jac_f^\dagger (p+\nu(\x-p))\|_\infty < n^2 \|D\Jac_f^\dagger(p)\|_\infty + {\rm O}(R).
\end{eqnarray}
We explicitly express $D \Jac_f^\dagger(\x)$ by computing the partial derivatives 
of each component of $\Jac_f^\dagger(\x)$. That is, 
\begin{eqnarray*}
D \Jac_f^\dagger(\x) = \frac{1}{\|\Jac_f(\x)\|_2^4} \left( 
\|\Jac_f(\x)\|_2^2 H_f(\x) - 2 \Jac_f(\x)^t \Jac_f(\x) H_f(\x)\right).
\end{eqnarray*}
We now  upper bound  $\|D \Jac_f^\dagger(p)\|_\infty$ by
\label{noncapisco2}
\begin{eqnarray*}\label{noncapisco2} 
\|D \Jac_f^\dagger(p)\|_\infty 
&\le& \frac{1}{\|\Jac_f(p)\|_2^4} \left( 
\|\Jac_f(p)\|_2^2  + 2 \|\Jac_f(p)^t\|_\infty \|\Jac_f(p)\|_\infty  \right)\|H_f(p)\|_\infty  \nonumber \\
&=& \frac{1}{\|\Jac_f(p)\|_2^4} \left( 
\|\Jac_f(p)\|_2^2  + 2 \|\Jac_f(p)^t\|_\infty \|\Jac_f(p)^t\|_1  \right)\|H_f(p)\|_\infty, \nonumber 
\end{eqnarray*}
where the last equality follows from Remark \ref{65}.
Proposition \ref{propMatEquivalence} 
 gives  $\|\Jac_f(p)^t\|_\infty \leq \|\Jac_f(p)^t\|_2$, as well as $\|\Jac_f(p)^t\|_1 \leq \sqrt{n} \|\Jac_f(p)^t\|_2$. Thus,
 noting that $\|\Jac_f(p)\|_2=\|\Jac_f(p)^t\|_2$, and using again Proposition \ref{propMatEquivalence}(2), one has
\begin{eqnarray}\label{stilettata} 
\|D \Jac_f^\dagger(p)\|_\infty
&\le& \frac{1}{\|\Jac_f(p)\|_2^4} \left( 
\|\Jac_f(p)\|_2^2  + 2 \sqrt{n} \|\Jac_f(p)\|^2_2  \right)\|H_f(p)\|_\infty \nonumber \\ 
&=& (1+2 \sqrt{n}) \frac{\|H_f(p)\|_\infty}{\|\Jac_f(p)\|_2^2}    
\le  (1+2\sqrt{n}) \frac{\|H_f(p)\|_\infty}{\|\Jac_f(p)\|_1^2}.  
\end{eqnarray}
By combining (\ref{supDJac}) and (\ref{stilettata}) and recalling that 
$\x \in \mathcal D$ implies $\|(\x-p)^t\|_\infty<R$, inequality (\ref{princDisugJac}) yields
$$
\|\Jac_f^\dagger(\x)\|_\infty <    
 \|\Jac_f^\dagger(p)\|_\infty  + n^2 (1+2\sqrt{n}) \frac{\|H_f(p)\|_\infty}{\|\Jac_f(p)\|_1^2}R+
{\rm O}(R^2).
$$
By definition of $\mathsf J$  we are done.
\qed

Next, we need a technical result (valid up to a first-order  error analysis).

\begin{lemma}\label{Full2} 
Let $f=f(\x)$ be a degree $\geq 2$ polynomial of~$P$ and 
let $p$ be a point of $\R^n$ such  that both  the Jacobian  $\Jac_f(p)$  and the Hessian matrix  $H_f(p)$  are nontrivial.
Let $R$ be a positive real number such that $R < \epsbold_{{\rm min}}$.
If  $R<\frac{\|\Jac_f(p)\|_1}{n^2 \|H_f(p)\|_\infty}$,
then, neglecting contributions of order ${\rm O}(R^2)$,
the Jacobian ${\rm Jac}_f(\x)$ is nonzero   for each  $\x
\in \sD(p,R)$.
\end{lemma}
\proof 
Since $R < \epsbold_{{\rm min}}$, we have $ \sD:=\sD(p,R)\subseteq {\bf B}(p)$.
  Lemma \ref{UppBoundH} applied  to the Hessian matrix $M(\x)=H_f(\x)$ gives, For each $\x\in \sD$,
$$
\sup_{0<\nu<1} \|H_f(p+\nu(\x-p))\|_\infty <n^2 \| H_f(p)\|_{\infty} + {\rm O}(R).
$$
On the other hand, inequality (\ref{depressione})  still holds true.
Thus we find
$
\|\Jac_f(\x)^t - \Jac_f(p)^t\|_\infty <n^2 \| H_f(p)\|_{\infty}R + {\rm O}(R^2)$,
so that  the same conclusion as in the  proof of Lemma \ref{Full1} gives now
$$
\|\Jac_f(\x)\|_1 > \|\Jac_f(p)\|_1 - n^2\|H_f(p)\|_\infty R -{\rm O}(R^2)>- {\rm O}(R^2).
$$ 
Therefore,  up to a first-order error analysis,  $\Jac_f(\bm x)$ 
is nonzero   for each $\x \in \mathcal D$, as we want.
\qed

As in the case  of necessary conditions,  it is natural to ask for a statement analogous  to Proposition \ref{PropPassa} which  avoids the computation of   the quantities ${\mathsf J}$ and ${\mathsf H}$,   providing  a  crossing cell condition simply  evaluating  
    the Jacobian and the Hessian matrix at a  given point.  The following statement  settles the question,   for small values of  the tolerance vector $\epsbold$, up to a  first-order error analysis, that is,  disregarding  second-order contributions. To this purpose, for the rest of this section, we assume  $\epsbold_{\rm max}\ll 1$.

\begin{prop}\label{PropPassa2}
Let $f=f(\x)$ be a degree $\geq 2$ polynomial of~$P$,
let $p$ be a point of $\R^n$ such that  Jacobian $\Jac_f(p)$ and the Hessian matrix $H_f(p)$ are nontrivial, 
and let ${\bf C}(p) \subseteq {\bf B}(p)$
be an $(\infty,\epsbold)$-unit cell  centered at $p$.
Let $R$ be a positive real number such that
$R < \min \left \{\epsbold_{\rm min}, \frac{\|\Jac_f(p)\|_1}{n^2 \|H_f(p)\|_\infty} \right\}$, 
let $c:=\max\{2, \sqrt{n}\}$ and set 
$$
\Theta:=  \|\Jac_f^\dagger(p)\|_\infty+  n^2(1+2 \sqrt{n}) \frac{\|H_f(p)\|_\infty}{\|\Jac_f(p)\|_1^2} R.
$$
If 
\begin{equation}
|f(p)| < \frac{2R}{\Theta (c+  n^{5/2} \|H_f(p)\|_\infty \Theta R)} =:B_2',
\end{equation}
then the hypersurface of equation $f=0$ crosses ${\bf C}(p)$  neglecting order ${\rm O}(R^2)$ contributions.
\end{prop}
\proof 
The proof runs parallel to that of Proposition \ref{PropPassa}. First, note that if $|f(p)|={\rm O}(R^2)$ there is nothing to prove, so we can assume 
$|f(p)|>{\rm O}(R^2)$.
From Lemma \ref{Full2} we know that  the Jacobian $\Jac_f(\bm x)$ is nonzero, 
for  $\x \in \mathcal D:=\mathcal D(p, R)$, up to a first-order analysis.
Moreover, since $R< \epsbold_{\rm min}$, one has $\sD\subseteq \mathcal {\bf C}(p)$.

We now consider the sequence of points $\{p_k\}$ and the column vectors $s_k$, $k\in \N$, defined  as in (\ref{DEF1}). 
Obviously, $p=:p_0\in \sD$. We prove by induction that,   for each $k\geq 1$, the points  $p_k$'s all lie in $ \mathcal D$  up to ${\rm O}(R^2)$, that is,
\begin{equation}\label{civuoleR4} 
\| (p_k-p)^t\|_\infty<R+{\rm O}(R^2),
\end{equation}
and satisfy the  inequality 
 \begin{equation}\label{civuoleR2}
 |f(p_k)| <  |f(p_{k-1})|+{\rm O}(R^2).\end{equation}

\noindent{{\bf Step I}} (The $k=1$ case). 
From the definitions of $s_0$ and $\Theta$ we have
$$
\|s_0\|_\infty = \|{\rm Jac}_f^\dagger(p)\|_\infty |f(p)| < \Theta |f(p)|.
$$
Furthermore,  since $c \geq 2$,  one  has by assumption
$|f(p)| <B_2' < \frac{2 R}{c \Theta} \le \frac{R}{\Theta}$.  Thus $\|s_0\|_\infty < R$, that is $p_1\in \sD$, showing condition (\ref{civuoleR4}) for $k=1$.

From the formal Taylor expansion of $f(\x)$ at $p$ (see 
Definition~\ref{defTaylor}) it follows that
$$
f(\x) = f(p) + \Jac_f(p) (\x-p)^t + \frac{1}{2}(\x-p)H_f(p) (\x-p)^t + R_2(\x),
$$
where $R_2(\bm x)$ is the $2$-nd remainder term of~$f(\x)$ at~$p$.
If we evaluate the former expression at~$p_1$ and apply 
Taylor's theorem (see Proposition~\ref{propTaylor}(1)) we get
\begin{eqnarray*}
f(p_1) = f(p) + \Jac_f(p) (p_1-p)^t + 
\frac{1}{2}(p_1-p)H_f(p) (p_1-p)^t + {\rm O}(R^3).
\end{eqnarray*}
Therefore, by  definitions (\ref{DEF1}) with $k=0$,
we get
\begin{eqnarray}\label{caldo2} \nonumber
f(p_1)  &=&  f(p) + \Jac_f(p) s_0 +  \frac{1}{2} s_0^t H_f(p) s_0 + {\rm O}(R^3)\\ \nonumber
&=&  f(p) - \frac{f(p)}{\|{\rm Jac}_f(p)\|_2^2} \Jac_f(p)  \Jac_f(p)^t +  
\frac{1}{2} \frac{|f(p)|^2}{\|{\rm Jac}_f(p)\|_2^4} \Jac_f(p) H_f(p) \Jac_f(p)^t + {\rm O}(R^3)\\ \nonumber
&=& f(p) - f(p) +  
\frac{1}{2} \frac{|f(p)|^2}{\|{\rm Jac}_f(p)\|_2^4} \Jac_f(p) H_f(p) \Jac_f(p)^t + {\rm O}(R^3)\\ 
&=&  |f(p)| \left ( 
\frac{1}{2} \frac{|f(p)|}{\|{\rm Jac}_f(p)\|_2^4} \Jac_f(p) H_f(p) \Jac_f(p)^t \right ) + {\rm O}(R^3). 
\end{eqnarray}
Let us upper bound the absolute value of the quantity
$$
Q:= \frac{1}{2} \frac{|f(p)|}{\|{\rm Jac}_f(p)\|_2^4} \Jac_f(p) H_f(p) \Jac_f(p)^t.
$$
To this end, use  H\"older's inequality, 
 Proposition  \ref{propMatEquivalence}(1) and Remark \ref{65}, to get:
\begin{eqnarray*}
|Q| &\le & \frac{1}{2} \frac{|f(p)|}{\|{\rm Jac}_f(p)\|_2^4}  \| \Jac_f(p)\|_2^2 \: \|H_f(p)\|_2\\
&\le & \frac{1}{2} \frac{|f(p)|}{\|{\rm Jac}_f(p)\|_2^2}  \sqrt{n} \|H_f(p)\|_\infty \\
&=&  \frac{1}{2} \sqrt{n} |f(p)| \frac{\|\Jac_f(p)^t\|_\infty}{\|{\rm Jac}_f(p)\|_2^2}  
\frac{\|H_f(p)\|_\infty}{\|\Jac_f(p)^t\|_\infty}\\
&=& \frac{1}{2} \sqrt{n} |f(p)| \; \|\Jac_f^\dagger(p)\|_\infty  \frac{\| H_f(p)\|_\infty}{\|\Jac_f(p)\|_1}.   
\end{eqnarray*}
Recalling the definition of $\Theta$ and the assumption on $R$, we find
$$
|Q|< \frac{1}{2} \sqrt{n} |f(p)| \Theta  \frac{ \| H_f(p)\|_\infty}{\|\Jac_f(p)\|_1} < \frac{1}{2} \sqrt{n} |f(p)| \frac{\Theta}{R}. 
$$

On the other hand, since $c\geq \sqrt{n}$,  one has by assumption  that 
$|f(p)| < B_2' < \frac{2R}{c \Theta} \le \frac{2R}{ \sqrt{n} \Theta}$.
Therefore  $|Q|<1$,  so that equality (\ref{caldo2})  proves condition (\ref{civuoleR2}) for $k=1$.
\smallskip

\noindent{{\bf Step II}} (Inductive step). Suppose that the points $p, p_1,\ldots,p_k$  of the sequence meet conditions (\ref{civuoleR4}) and (\ref{civuoleR2}). Further assume  $|f(p_i)|>{\rm O}(R^2)$, whence, in particular, $\|s_i\|_\infty=\|\Jac_f^\dagger (p_i)\|_\infty|f(p_i)|>{\rm O}(R^2)$, $i=1,\ldots,k$. We want to show that the point $p_{k+1}$ satisfies such conditions as well.
From the formal Taylor expansion of $f(\x)$ at $p_{i-1}$ (see 
Definition~\ref{defTaylor}) it follows that
$$
f(\x) = f(p_{i-1}) + \Jac_f(p_{i-1}) (\x-p_{i-1})^t + \frac{1}{2}(\x-p_{i-1})H_f(p_{i-1}) (\x-p_{i-1})^t + R_2(\x),
$$
where $R_2(\bm x)$ is the $2$-nd remainder term of~$f(\x)$ at~$p_{i-1}$, $i=1,\ldots,k$.
If we evaluate the former expression at~$p_i$ and apply 
Taylor's theorem (see Proposition~\ref{propTaylor}(1)) we get
\begin{equation}\label{EqT1bis}
f(p_i) = f(p_{i-1}) + \Jac_f(p_{i-1}) (p_i-p_{i-1})^t + 
\frac{1}{2}(p_i-p_{i-1})H_f(p_{i-1}) (p_i-p_{i-1})^t + {\rm O}(R^3).
\end{equation}
Then  the same argument  giving inequality (\ref{idem}) in the proof of Proposition \ref{PropPassa}  yields now
$$
|f(p_i)|\leq  \frac{1}{2}  \|H_f(p_{i-1})\|_1 \|(p_i - p_{i-1})^t\|_\infty^2 + {\rm O}(R^3).
$$
Since the Hessian matrix is symmetric one has
$\|H_f(p_{i-1})\|_1=\|H_f(p_{i-1})\|_\infty$. Moreover, from  Lemma \ref{UppBoundH}
applied to  $M(\x)= H_f(\x)$, we have $\|H_f(p_{i-1})\|_\infty < n^2 \|H_f(p)\|_\infty + {\rm O}(R)$. 
Thus the previous relation yields
\begin{eqnarray}\label{ineq1bis}
|f(p_i)| &<&   
\frac{1}{2} n^2 \|H_f(p)\|_\infty \|s_{i-1}\|_\infty^2 + \frac{1}{2}  \|s_{i-1}\|_\infty^2 {\rm O}(R)+ {\rm O}(R^3) \nonumber\\
&\le& \frac{1}{2} n^2 \|H_f(p)\|_\infty \|s_{i-1}\|_\infty^2 + {\rm O}(R^3),
\end{eqnarray}
where 
the last inequality is a consequence of  the inductive assumption $$\| (p_i-p_{i-1})^t\|_{\infty}\leq \| (p_i-p)^t\|_{\infty}+\| (p_{i-1}-p)^t\|_{\infty}<2R+{\rm O}(R^2), \;\; i=1,\ldots, k.$$ Proposition \ref{1Luglio} reads ${\mathsf J}<\Theta+{\rm O}(R^2)$. Then 
inequality (\ref{ineq1bis})  gives 
\begin{eqnarray*}
\|s_i\|_\infty &=& \|\Jac_f^\dagger(p_i)\|_\infty |f(p_i)|
< (\Theta +{\rm O}(R^2)) \left( \frac{1}{2} n^2 \|H_f(p)\|_\infty \|s_{i-1}\|_\infty^2 + {\rm O}(R^3) \right)\\
&\le& \frac{1}{2} n^2 \Theta \|H_f(p)\|_\infty  \|s_{i-1}\|_\infty^2 + {\rm O}(R^3).
\end{eqnarray*}
Now, define $\tau_i : =\frac{\|s_i\|_\infty}{\|s_{i-1}\|_\infty}$.  Recalling  the inductive assumption $ |f(p_i)| <  |f(p_{i-1})|+{\rm O}(R^2)$, $i=1,\ldots,k$,
it thus follows that 
\begin{eqnarray}\label{ineqTau}
\tau_i &=& \frac{\|s_i\|_\infty}{\|s_{i-1}\|_\infty} 
\le \frac{1}{2} n^2 \Theta \|H_f(p)\|_\infty  \|s_{i-1}\|_\infty + \frac{{\rm O}(R^3)}{\|s_{i-1}\|_\infty} \nonumber\\
&<& \frac{1}{2} n^2 \Theta \|H_f(p)\|_\infty  \|s_{i-1}\|_\infty + {\rm O}(R) \nonumber \\
&<& \frac{1}{2} n^2 \Theta^2 \|H_f(p)\|_\infty  |f(p_{i-1})| + {\rm O}(R) 
< \frac{1}{2} n^2 \Theta^2 \|H_f(p)\|_\infty  |f(p)| +{\rm O}(R). 
\end{eqnarray}

Setting  $T:=\frac{1}{2} n^2 \Theta^2 \|H_f(p)\|_\infty  |f(p)|$, and recalling that
 $$|f(p)| < B_2' < \frac{2}{n^{5/2} \|H_f(p)\|_\infty \Theta^2 }, $$ we then have
$T <\frac{1}{\sqrt{n}}\leq 1$,
so that  $\tau_i<T + {\rm O}(R)$ by  inequality (\ref{ineqTau}).
Thus, as in  the corresponding step in the proof of Proposition \ref{PropPassa}, we get the bound:
\begin{eqnarray*}
\|(p_{k+1} -p)^t\|_\infty 
&\leq & \|s_0\|_\infty (1 + \tau_1+\tau_1\tau_2 +\ldots + \tau_1\tau_2 \dots \tau_k)\\
&<& \|s_0\|_\infty \sum_{i=0}^k (T + {\rm O}(R))^i = 
\|s_0\|_\infty \left (\sum_{i=0}^k T^i +  {\rm O}(R)\right)\\
&\leq& \|s_0\|_\infty  \sum_{i=0}^\infty T^i + {\rm O}(R^2)  =\frac{\|s_0\|_\infty}{1- T} + {\rm O}(R^2)\\
&=& \frac{\| \Jac_f^\dagger(p)\|_\infty |f(p)|}{1- T} + {\rm O}(R^2)<
 \frac{\Theta |f(p)|}{1- T} + {\rm O}(R^2).
\end{eqnarray*}
By  definition of $T$ 
we then find
\begin{eqnarray*}
\|(p_{k+1} -p)^t\|_\infty &<& 
\frac{\Theta |f(p)|}{1- \frac{1}{2} n^2 \Theta^2 \|H_f(p)\|_\infty |f(p)|} +{\rm O}(R^2)\\ 
&=& \frac{2 \Theta |f(p)|}{2- n^2 \Theta^2 \|H_f(p)\|_\infty |f(p)|} +{\rm O}(R^2)
< R + {\rm O}(R^2),
\end{eqnarray*}
where the last inequality rewrites as $|f(p)|<\frac{2R}{\Theta(2+n^2\Theta \|H_f(p)\|_\infty R)}$, which follows  from the assumption $|f(p)|<B_2'$.
We then conclude that $p_{k+1} \in \mathcal D$ (up to ${\rm O}(R^2$)). 

Moreover,   
we observe that relation (\ref{caldo2}) can be easily adapted to the pair of points $p_k$, $p_{k+1}$
in the form
\begin{eqnarray}\label{911}
f(p_{k+1}) = |f(p_k)| \left ( 
\frac{1}{2} \frac{|f(p_k)|}{\|{\rm Jac}_f(p_k)\|_2^4} \Jac_f(p_k) H_f(p_k) \Jac_f(p_k)^t \right ) + {\rm O}(R^3). 
\end{eqnarray}
Let us upper bound the absolute value of the quantity
$$
Q_k:= \frac{1}{2} \frac{|f(p_k)|}{\|{\rm Jac}_f(p_k)\|_2^4} \Jac_f(p_k) H_f(p_k) \Jac_f(p_k)^t.
$$
As in Step I, by using H\"older's inequality
and Proposition \ref{propMatEquivalence}(1), we obtain
$$
|Q_k| 
\le  \frac{1}{2} \sqrt{n} \frac{|f(p_k)|}{\|{\rm Jac}_f(p_k)\|_2^2}   \|H_f(p_k)\|_\infty 
\le  \frac{1}{2}\sqrt{n}  \frac{|f(p_k)|}{\|{\rm Jac}_f(p_k)\|_1^2}  \ \|H_f(p_k)\|_\infty.
$$
By Lemma \ref{UppBoundH} applied to the matrix $M(\x) = H_f(\x)$, we then have
\begin{eqnarray}\label{Qk}
|Q_k| 
&< & \frac{1}{2}\sqrt{n}  \frac{|f(p_k)|}{\|{\rm Jac}_f(p_k)\|_1^2}  (n^2\|H_f(p)\|_\infty+{\rm O}(R)) \nonumber\\
&\le& \frac{1}{2}n^{5/2}  \frac{|f(p_k)|}{\|{\rm Jac}_f(p_k)\|_1^2}   \|H_f(p)\|_\infty+{\rm O}(R^2), 
\end{eqnarray}
where the last inequality comes from $|f(p_k)| < |f(p)| + {\rm O}(R^2) < \frac{R}{\Theta} + {\rm O}(R^2)$, 
a consequence of the inductive hypothesis and the assumption $|f(p)|<B_2'$.
By H\"older's inequality and the definition of $\Theta$ we find
$$
1= |\Jac_f(p_k) \Jac_f(p_k)^\dagger| \le \|\Jac_f(p_k)\|_1  \|\Jac_f(p_k)^\dagger\|_\infty
< \|\Jac_f(p_k)\|_1(\Theta + {\rm O}(R^2)). 
$$ 
Therefore inequality (\ref{Qk}) becomes:
\begin{eqnarray*}
|Q_k| &< & \frac{1}{2} n^{5/2} |f(p_k)| \|H_f(p)\|_\infty(\Theta+{\rm O}(R^2))^2+{\rm O}(R^2)\\
&<&  \frac{1}{2} n^{5/2} ( |f(p)| + {\rm O}(R^2)) \|H_f(p)\|_\infty \Theta^2+{\rm O}(R^2)\\
&=& \frac{1}{2} n^{5/2}  \|H_f(p)\|_\infty \Theta^2  |f(p)|+{\rm O}(R^2).
\end{eqnarray*}
On the other hand,  
$|f(p)| < B_2' < \frac{2}{n^{5/2} \|H_f(p)\|_\infty  \Theta^2}$.
Thus we  find $|Q_k|<1 + {\rm O}(R^2)$,  so that equality (\ref{911})  
yields  the desired  condition $|f(p_{k+1})| < |f(p_k)|+{\rm O}(R^2)$.
\smallskip

\noindent{{\bf Step  III}} (Conclusion). If there exists $k\in \N$ such that $|f(p_k)|={\rm O}(R^2)$ we are done. Otherwise,
 we know from Step II that $\tau_k:=\frac{ \|s_k\|_\infty}{ \|s_{k-1}\|_\infty}<T +{\rm O}(R)<1+{\rm O}(R)$ for $k\in \N$, whence $\tau_k<1$ for $R\ll1$.  Then, the same argument as in Step III of the proof of Proposition \ref{PropPassa} applies to say that
the sequence of points $\{p_k\}_{k \in \N}$ converges to a point $p^*$. Since the  $p_k$'s belong to $\sD$ up to ${\rm O}(R^2)$, the point $p^*$ belongs to the closure $\overline{\sD}\subseteq {\bf C}(p)$ up to ${\rm O}(R^2)$.
We also know that $\| s_k\|_\infty=\tau_1\tau_2\ldots\tau_k\|s_0\|_\infty<\tau^k\|s_0\|_\infty$, where $\tau=\sup_{k\in \N}\{\tau_k\}$. Therefore $\lim_{k\to \infty}\|s_k\|_\infty< \lim_{k\to \infty}\tau^k\|s_0\|_\infty=0$.
From  inequality~(\ref{ineq1bis}), 
 we then conclude that 
$$|f(p^*)|=\lim_{i\to \infty} |f(p_i)|\le\frac{1}{2} n^2 \|H_f(p)\|_\infty \lim_{i\to\infty}\|s_{i-1}\|_\infty^2 + {\rm O}(R^3)={\rm O}(R^3),$$  so that the hypersurface $f=0$ crosses the cell ${\bf C}(p)$  neglecting  order ${\rm O}( R^3)$ (hence order ${\rm O}(R^2))$ contributions.  \qed

\begin{example*}\label{EE}
We consider the polynomial $f(x,y)=y^2+x^3-x-3$ and the point 
$p=(0, 1.7)$.
We let $\epsbold =(0.06,0.06)$ and consider the
$(\infty,\epsbold)$-unit ball   ${\bf B}(p)$ centered at $p$.
Direct computations show that $\|\Jac_f(p)\|_1=3.4 \neq 0$ and
$\|H_f(p)\|_\infty=2$, whence $\frac{\|\Jac_f(p)\|_1}{n^2 \|H_f(p)\|_\infty }=0.425$.
 We choose $R=0.05 < \min\{0.06, 0.425\}$.
Further, we compute $\Theta \approx 0.403$, so for the bound $B_2'$ 
of Proposition \ref{PropPassa2}  we find  
 $$
 B_2' = \frac{2R}{\Theta (2+ n^{5/2}\|H_f(p)\|_\infty \Theta R)} \approx 0.1113.
$$
Since $|f(p)|=\frac{11}{10}=0.11$ is strictly smaller than $B_2'$, 
by using Proposition \ref{PropPassa2} we conclude that the elliptic
curve of equation $f=0$ crosses the  ball ${\bf B}(p)$. 
\end{example*}

Keeping the assumptions and notation as in propositions \ref{PropNonPassa2},  \ref{PropPassa} and \ref{PropPassa2}
we conclude this section 
comparing the bounds $B_1'$, $B_2$ and $B_2'$ provided by such propositions.

\begin{prop}\label{PropB1B2} Notation and assumptions as above. Further assume $\epsbold_{{\rm max}}\ll 1$ and let $R$ be a positive real number such that
$$
R < \min\left \{\epsbold_{{\rm min}}, \frac{\|\Jac_f(p)\|_1}{\mathsf H},
\frac{\|\Jac_f(p)\|_1}{n^2 \|H_f(p)\|_\infty} \right\}.
$$  Then 
$B_2' +{\rm O}(R^3) < B_2 < B_1' $.
\end{prop}
\proof 
Recalling the definitions of the quantities
$B_2$, $\mathsf J$, $\mathsf H$, $\Theta$  we get
$$
B_2 = \frac{2R}{{\mathsf J}(c+\sqrt{n}\mathsf H {\mathsf J} R)}  >
\frac{2R}{(\Theta + {\rm O}(R^2))\big(c+\sqrt{n}\mathsf H (\Theta + {\rm O}(R^2)) R\big)} 
= \frac{2R}{\Theta (c+\sqrt{n}\mathsf H \Theta R)+ {\rm O}(R^2)}. 
$$
  Lemma \ref{UppBoundH}, applied  to the Hessian matrix  in the open convex set $\sD=\sD(p,R)$, yields $\mathsf H < n^2 \|H_f(p)\|_\infty + {\rm O}(R)$, so that 
\begin{eqnarray}\label{lafine} \nonumber
B_2&>&   \frac{2R}{\Theta \Big(c+\sqrt{n}\big(n^2 \|H_f(p)\|_\infty+{\rm O}(R)\big) \Theta R\Big)+ {\rm O}(R^2)} \\  
&=& \frac{2R}{\Theta \big(c+n^{5/2} \|H_f(p)\|_\infty \Theta R\big)+ {\rm O}(R^2)}.
\end{eqnarray}	
Write the right-hand side term of (\ref{lafine}) as
\begin{equation}\label{evicina}
\frac{2R}{\Theta (c+n^{5/2} \|H_f(p)\|_\infty\Theta R)}\; \frac{1}{1+\frac{{\rm O}(R^2)}{c+n^{5/2} \|H_f(p)\|_\infty\Theta R}}.
\end{equation}
Noting that
$\left | \frac{{\rm O}(R^2)}{c+n^{5/2} \|H_f(p)\|_\infty\Theta R} \right |<1$ for $R  \ll 1$, the quantity in (\ref{evicina}) rewrites as
\begin{multline*} 
\frac{2R}{\Theta (c+n^{5/2} \|H_f(p)\|_\infty\Theta R)}\left(1-\frac{{\rm O}(R^2)}{c+n^{5/2} \|H_f(p)\|_\infty\Theta R}+\cdots\right)  \\
  = \frac{2R}{\Theta (c+n^{5/2} \|H_f(p)\|_\infty\Theta R)} + {\rm O}(R^3),
\end{multline*}
By  definition of $B_2'$, the lower bound of $B_2$ then easily follows.

From Lemma~\ref{Full2} we know that the Jacobian 
$\Jac_f(\bm x)$ is nonzero in $\mathcal D$,
so  that 
$\Jac_f(p) \Jac_f^\dagger(p) =1$. Thus  H\"older's inequality yields
$
1  \le \| \Jac_f(p)\|_1 \| \Jac_f^\dagger(p)\|_\infty \le
\|\Jac_f(p)\|_1\mathsf J 
$.
Therefore, recalling  the assumption on $R$,  we  upper bound  $B_2$  by
\begin{equation}\label{serve}
B_2 = \frac{2R}{\mathsf J(c+\sqrt{n}{\mathsf H} \mathsf J  R)} < \frac{2R}{c\mathsf J} 
\le \frac{R}{\mathsf J} \le \|\Jac_f(p)\|_1 R < \|\Jac_f(p)\|_1 \epsbold_{\rm max}.
\end{equation}
Recalling the definition of $B_1'$ we get the desired upper bound.
\qed

\begin{rem*}\label{REM} The inequality  $B_1 > B_2$ holds true. 
For degree $\geq 2$ polynomials, it is a consequence of  Lemma \ref{PropB1B1'} and Proposition \ref{PropB1B2}; in the linear case 
it follows  from Remark \ref{lineare} and equality (\ref{serve}).
It then may happen that $|f(p)|$ belongs to the interval $(B_2, B_1)$.
In this case, with the only use of propositions~\ref{PropNonPassa}
and~\ref{PropPassa}, we cannot conclude whether or not  the hypersurface $f=0$ crosses a given  unit cell ${\bf C}(p)$.
Because of the local nature of the previous results,
a more accurate  analysis, performed by iteratively considering smaller cells, may overcome 
that problem. Correspondingly, up to a second-order analysis, Proposition \ref{PropB1B2} gives  $B_1' > B_2'$.
\end{rem*}

\section{The crossing area algorithm}
\label{RAL}\addtocounter{subsection}{1}\setcounter{theorem}{0}

We keep the notation of the previous sections. In particular,  letting $\x=(x_1,\ldots,x_n)$ be indeterminates, we recall 
that $f=f(\x)$ is a non-costant polynomial of $P$, 
$p$ is a point of $\mathbb R^n$, $\varepsilon_1,\ldots,\varepsilon_n$
are positive real numbers, ${\bm \varepsilon}$ denotes the  tolerance   vector $(\varepsilon_1,\ldots,\varepsilon_n)$, 
and ${\bf C}(p) \subseteq {\bf B}(p)$ is an $(\infty,{\bf \varepsilon})$-unit cell centered
at $p$ (see Section \ref{passaggioCella1}). 
In this section, we describe  an algorithm to decide whether or not 
a hypersurface of equation $f=0$  intersect a given bounded region
in the affine space ${\mathbb A}_{\x}^n(\reals)$.

Let's start with a local version of  such a  crossing problem. That is, we introduce
an algorithm that, given a non-costant polynomial $f$ and a unit cell ${\bf C}(p)$, 
returns a value which describes the intersection of the hypersurface of equation $f=0$ 
with  ${\bf C}(p)$. Namely,
\begin{itemize}
\item $0$ if the hypersurface  $f=0$ does not cross ${\bf C}(p)$;
\item $1$ if the hypersurface  $f=0$ crosses ${\bf C}(p)$;
\item $\zeta$ (unknown) if neither Proposition \ref{PropNonPassa} nor Proposition  \ref{PropPassa} applies.
\end{itemize}

Summarizing, we have:

\medskip
\noindent
\framebox{
\parbox[l]{0.98\textwidth}
{{\bf The CROSSING CELL algorithm}\\
Given a non-costant polynomial $f=f(\x) \in P$, 
a point $p \in \mathbb R^n$ such that  both the  Jacobian  $\Jac_f(p)$  and the Hessian matrix  $H_f(p)$  are nontrivial at  $p$, and 
a  tolerance vector ${\bf \varepsilon} = (\varepsilon_1,\ldots,\varepsilon_n)$,
the algorithm returns an element of $\{0,1,\zeta\}$.
\begin{description}
\item[1.] Compute $|f(p)|$, and the bounds $B_1$ and $B_2$ from  propositions \ref{PropNonPassa} 
and \ref{PropPassa} (see also Remark \ref{lineare}).
\item[2.] If $|f(p)| > B_1$ return $0$;
if $|f(p)| < B_2$ return $1$;
else return $\zeta$.
\end{description}
}}

\begin{rem*}\label{REM2} We observe that there may be variants of the previous algorithm. For instance:
\begin{enumerate}
\item  The crossing cell algorithm could be performed up to a first-order error analysis
       simply replacing the bounds $B_1$ and $B_2$ by $B_1' $ and $ B_2'$ (defined in propositions \ref{PropNonPassa2} and \ref{PropPassa2}).
 \item   In the case $B_2<f(p)<B_1$, in order to limit the problem of indeterminacy (already pointed out in Remark \ref{REM}), 
      the crossing cell algorithm could be performed by iteratively considering
      smaller unit cells ${\bf C}(p)$, simply obtained by a subdivision procedure
     ending as soon as  such  cells are sufficiently small  to solve  the problem in the given context.
\end{enumerate}
\end{rem*}

Now, we consider a more general version of the crossing problem. 
Let $f=f(\x)\in P$  as above, and let $\mathcal T$ be a nontrivial bounded region 
of ${\mathbb A}_{\x}^n(\reals)$ of   type
$
\mathcal T := [a_1,b_1] \times \cdots \times [a_n,b_n],
$
with $a_k, b_k \in \reals$ and $a_k < b_k$, for each $k=1,\ldots,n$.
In order to reduce the problem of studying the intersection  of the hypersurface $f=0$  with the region $\mathcal T$
to a local crossing cell problem, we perform a discretization of 
the region $\mathcal T$ as follows.
Let $d_k$ be  the {\em sampling distance} with respect to the  component $x_k$, $k=1,\ldots,n$.
For each $k=1,\ldots,n$, we define 
\begin{equation}\label{Jkxk}
 J_k := \left \lceil \frac{b_k-a_k- \frac{d_k}{2}}{d_k} \right \rceil + 1 \; \;\;\; {\rm and} \;\;\; \;x_{k,j_k} :=a_k + j_k d_k,   
\end{equation}
where $\lceil x \rceil = \min\{z \in \mathbb N \:|\: z  \ge x\}$ and
$j_k=0,\ldots, J_k-1$. Here $J_k$ denotes the number of considered samples for each component, and $j_k$ the index of the sample.
We denote by ${\bf j}$ the multi-index $(j_1,\ldots, j_n)$,
by ${\x}_{\bf j}:=(x_{1,j_1},\ldots,x_{n, j_n})$ the ${\bf j}$-th {\em sampling point}, 
and by
\begin{equation*}\label{cells}
{\bf C}({\bf j}):=\left\{ {\x} =( x_1,\ldots, x_n ) \in \mathbb A^n(\reals) \; \Big| \;  
x_k\in \left[x_{k,j_k}-\frac{d_k}{2}, x_{k,j_k}+\frac{d_k}{2}\right),\; k=1,\ldots,n \right\} 
\end{equation*}
{\em the cell centered at (and represented by) the point} ${\x}_{\bf j}$.
The discretization of $\mathcal T$ is given by the $J_1\times \cdots \times J_n$ cells 
of type ${\bf C}({\bf j})$ which are  a covering of the region $\sT$.
Let us stress the fact that the discretization is defined by relations (\ref{Jkxk}), that is, 
by the initializing point $a=(a_1,\ldots,a_n) \in \sT$ and the discretization step 
$d=(d_1,\ldots,d_n)$. 

 We define 
a multi-matrix $\sA=(\gra_{j_1j_2\cdots j_n})$ of type $ J_1\times \cdots \times J_n$, 
called the {\em crossing area  matrix} (with respect to the given discretization $\{a,d\}$ defined by $a$ and $d$).   
Each entry $\sA({\bf j})$ of $\sA$ contains informations about  whether or not  the hypersurface  of equation  $f=0$  crosses the cell ${\bf C}({\bf j})$.
Namely,
\begin{itemize}
\item $\sA({\bf j})=0$ if the hypersurface $f=0$ does not cross ${\bf C}({\bf j})$;
\item $\sA({\bf j})=1$ if the hypersurface  $f=0$ crosses ${\bf C}({\bf j})$;
\item $\sA({\bf j})=\zeta$ (unknown) if neither Proposition \ref{PropNonPassa} nor Proposition  \ref{PropPassa} applies.
\end{itemize}

In order to use propositions \ref{PropNonPassa} and \ref{PropPassa},  
we interpret each cell ${\bf C}({\bf j})$ in terms of the unit ball of an appropriate 
normed space $\reals^n$, where the appropriate norm needs to be defined. 
To this end, let 
$\epsbold=(\varepsilon_1,\ldots, \varepsilon_n):=\big(\frac{d_1}{2},\ldots,\frac{d_n}{2} \big)=:\frac{d}{2}$,
let ${\mathcal E} \in \textrm{Mat}_n(\reals)$ be the positive diagonal matrix with entries
$1/\varepsilon_1, \ldots, 1/\varepsilon_n$, and consider the ${\mathcal E}$-weighted $\infty$-norm 
on $\reals^n$ (see Definition \ref{weightedNorm}).
Since
$||({\x}-{\x}_{\bf j})^t||_{\infty,\sE}= \| \sE({\x}-{\x}_{\bf j})^t\|_{\infty}=\max_{k=1,\ldots,n}\Big\{\big|\frac{2(x_k-x_{k,j_k})}{d_k}\big|\Big\}$,
we can express the $(\infty,\varepsilon )$-unit ball centered at $\x_{\bf j}$  as  $${\bf B}({\bf j})=\Big\{ {\x}=(x_1,\ldots,x_n)\in {\mathbb A}^n(\reals)\; \big|\;  \max_{k=1,\ldots,n} \big\{|(x_k-x_{k,j_k}|\big\}\leq \frac{d_k}{2} \Big\}.$$  
Therefore
$$
{\bf C}({\bf j})={\bf B}({\bf j})\setminus \bigcup_{k=1,\ldots,n}L^{+}_k,
$$ 
where $L^{+}_k$ is the hyperplane of equation 
$L^{+}_k: x_k = x_{k,n_k}+\varepsilon_k, \;\;k =1,\ldots,n$. 
According to the notation settled at the beginning of Section \ref{passaggioCella1}, 
it follows that ${\bf C}({\bf j})$ is an $(\infty,\varepsilon )$-unit cell centered at ${\x}_{\bf j}$.
The whole procedure is gathered in the following algorithm.

\medskip
\noindent
\framebox{
\parbox[l]{0.98\textwidth}
{{\bf The CROSSING AREA algorithm}\\
Given a non-costant polynomial $f=f(\x) \in P$, 
a region $\mathcal T := [a_1,b_1] \times \cdots \times [a_n,b_n]$
of ${\mathbb A}_{\x}^n(\reals)$, and 
a discretization step $d = (d_1,\ldots,d_n)$,
the algorithm returns a multi-matrix $\sA$ with values in $\{0,1,\zeta\}$.

\begin{description}
\item[1.] Let $a=(a_1,\ldots,a_n) \in \mathcal T$ be the initializing point. 
Using relations (\ref{Jkxk})  construct  the  discretization $\{a,d\}$
of the region $\mathcal T$ and the multi-matrix $\sA$ of size 
$J_1 \times \cdots \times J_n$.

\item[2.] For each multi-index ${\bf j}=(j_1,\ldots,j_n)$
assign to $\sA({\bf j})$ the output of the CROSSING CELL algorithm
applied to the polynomial $f$, the point $p=x_{\bf j}$ and the tolerance vector  $\varepsilon=\frac{d}{2}$.

\item[3.] Return $\sA$. 
\end{description}
}}

\medskip
An implementation of the  CROSSING CELL and the CROSSING AREA algorithms
has been done using \cocoa 5 (see \cite{Cocoa}) and is available at 
{\tt http://www.dima.unige.it/$\sim$torrente/recognitionAlgorithm.cocoa5.}

\section{An application to the  Hough transform}\label{RecAl}
\addtocounter{subsection}{1}\setcounter{theorem}{0}

In this section we discuss an explicit manner of how using 
the CROSSING AREA algorithm in a case of special interest, 
that is,  the Hough transform technique. The Hough transform is a pattern recognition technique,
based on algebraic geometry arguments, for the automated recognition of curves in images.
We refer to  \cite{DH},  \cite{BMP} and  \cite{BR}  for  background material and complete details, 
and to \cite{etal}, \cite{Ro} for applications and further developments. Here, we restrict ourselves 
to just recall few basic definitions and properties.

Most of the results in this section  hold over an infinite integral ring $K$. 
However, let us restrict to the case of interest in the applications, assuming $K= \reals$.

For every $t$-tuples of independent  parameters ${\bm \lambda}:=(\lambda_1,\ldots,\lambda_t)\in \R^t$, let 
\begin{equation}\label{GE} 
f_{\bm \lambda}({\x})=\sum_{i_1,\ldots,i_n} x_1^{i_1} \cdots  x_n^{i_n}g_{i_1\cdots i_n}({\bm \lambda}) , \;\;\;\; i_1+\cdots+i_n\leq d,
\end{equation}
be a family ${\mathcal P}$ of  irreducible polynomials in the indeterminates ${\x}:=(x_1,\ldots,x_n)$,  
of a given degree $d$ (not depending on ${\bm \lambda}$), whose coefficients $g_{ij}({\bm \lambda})$ 
are  expressed polynomially in ${\bm \lambda}$. Let $\mathcal{F}$ be the corresponding  family of the zero loci 
$\mathcal{H}_{\bm \lambda} = \{\x \in {\mathbb A}_{\x}^n(\R) \:|\: f_{\bm \lambda}(\x)=0\} $, 
and let assume that $\mathcal{H}_{\bm \lambda}$ is a hypersurface for each parameter ${\bm \lambda}$ belonging to an euclidean open subset $\sU\subseteq \R^t$
(of course, this is always  the  case  if  the base field $K$ is  algebraically closed). 
Clearly, such hypersurfaces are irreducible, since  the polynomials of the family ${\mathcal P}$ 
are assumed to be irreducible in $\R[{\x}]$. So, we want  $\mathcal{F}$ to be a  
{\em family of irreducible hypersurfaces which share the degree}.

\begin{definition*}\label{def1}
Let  $\mathcal{F}$ be a family of hypersurfaces  ${\mathcal{H}}_{\bm \lambda}$ as above, 
and let $p=(x_1(p),\ldots, x_n(p))$ be a point in the {\em image space}  ${\mathbb A}_{\x}^n(\R)$. 
Let ${\bm \Lambda}:=(\Lambda_1,\ldots,\Lambda_t)$ be indeterminates, and let
${\Gamma}_p(\mathcal{F})$ be the hypersurface defined in the affine $t$-dimensional {\em parameter space} 
${\mathbb A}_{\bm \Lambda}^t(\R)$ by the polynomial equation 
$$f_p(\bm \Lambda) =
\sum_{i_1,\ldots,i_n} x_1(p)^{i_1} \cdots  x_n(p)^{i_n}g_{i_1\cdots i_n}({\bm \Lambda})=0, \;\;\;\;\;\; i_1+\cdots+i_n\leq d. 
$$ 
We say that ${\Gamma}_p(\mathcal{F})$ is the {\em Hough transform of the point $p$ 
with respect to the family} $\mathcal{F}$.  
If no confusion will arise, we simply say that ${\Gamma}_p(\mathcal{F})$ 
is the {\em Hough transform of $p$}. 
\end{definition*}
\smallskip

Summarizing,  the polynomials family defined by (\ref{GE})  gives rise to a polynomial 
$F({\x}; {\bm \Lambda})\in \R[{\x}; {\bm \Lambda}]$ giving, for each ${\bm \lambda} \in\sU$ 
and for each point $p \in {\mathbb A}_{\x}^n(\R)$, 
$$
{\mathcal{H}}_{\bm \lambda} : F({\x};{\bm \lambda})=f_{\bm \lambda}({\x})=0 \;\;\;\;
{\rm and}
\;\;\;\;
{\Gamma}_p(\mathcal{F}) : F(p; {\bm \Lambda})=f_p(\bm \Lambda)=0. 
$$

The following general facts hold true (see \cite[Theorem 2.2, Lemma 2.3]{BMP}, \cite[Section 3]{BR}).
 \begin{enumerate}
\item[1.]
{\em The Hough transforms ${\Gamma}_p(\mathcal F)$ of the pairs $(\mathcal{H}_{\bm \lambda}, p)$, 
when $p$ varies on $\mathcal{H}_{\bm \lambda}$, all pass through the point ${\bm \lambda}$}.
\item[2.] {\em Assume that  the Hough transforms ${\Gamma}_p(\mathcal F)$, 
when $p$ varies on $\mathcal{H}_{\bm \lambda}$,  have  a point in common other than ${\bm \lambda}$,
say ${\bm \lambda^{\prime}}$. Thus the two hypersurfaces $\mathcal{H}_{\bm \lambda}$, 
$\mathcal{H}_{\bm \lambda^{\prime}}$ coincide}.
\item[3.]  (Regularity property)  {\em The following conditions are equivalent}:
\begin{enumerate}
\item[(a)] {\em For  any  hypersufaces  $\mathcal{H}_{\bm \lambda}$,  $\mathcal{H}_{\bm \lambda^{\prime}}$ 
in $\mathcal{F}$, the equality  $\mathcal{H}_{\bm \lambda}=\mathcal{H}_{\bm \lambda^{\prime}}$ 
implies ${\bm \lambda}={\bm \lambda^{\prime}}$}.
\item[(b)] {\em For each  hypersurface $ \mathcal{H}_{\bm \lambda}$ in $\mathcal{F}$, one has
$\bigcap_{p\in \mathcal{H}_{\bm \lambda}}{\Gamma}_p(\mathcal F)={\bm \lambda}$}.
\end{enumerate} 
\end{enumerate}

A family $\mathcal{F}$ which meets one of the above equivalent conditions is said to be
 {\em Hough regular}.

 \smallskip

Let us consider the case $n=2$ we are interested in. Given a profile of interest in the image plane ${\mathbb A}_{\x}^2(\R)$, 
the Hough approach detects a curve of the family $\mathcal F=\{\mathcal{H}_{\bm \lambda}\}$ 
 best approximating it. 
Under the assumption of Hough regularity on $\mathcal F$, the  
detection procedure can be  highlighted as follows. 
\begin{enumerate}
\item[I.] Choose  a set $\mathbb X = \{p_1,\ldots,p_\nu\}$ of points of interest  in the image plane $ {\mathbb A}_{\x}^2(\reals)$.
\item[II.] In the parameter space ${\mathbb A}_{\bm \Lambda}^t(\R)$ 
find the (unique) intersection of the Hough transforms corresponding 
to the points $p_\iota$, that is, compute
${\bm \lambda} = \bigcap_{\iota=1,\ldots,\nu} \Gamma_{p_\iota}(\mathcal F)$. 
\item[III.] Return the  curve ${\mathcal H}_{\bm \lambda}$ uniquely determined by the parameter ${\bm \lambda}$.
\end{enumerate}

From a practical point of view, Step II is usually performed using the 
so called ``voting procedure", a discretization  approach for the (not easy) problem 
of computation of the intersection point ${\bm \lambda}$. 
Its core consists of the following three steps. Find a proper discretization of a suitable bounded 
region $\sT $ contained in the  open set  $\sU\subset \R^t$
 of the parameter space. Construct on it an {\em accumulator function}, that is,  
 a function that, for each Hough transform  $\Gamma_{p_\iota}(\mathcal F)$ and for each cell 
of the discretized region, records and sums the ``vote" $1$, if $\Gamma_{p_\iota}(\mathcal F)$ 
crosses the cell, and the ``vote" $0$ otherwise.   
Optimize the accumulator function by computing the cell corresponding to the (local) maximum; as suggested by the general results  
 recalled above, the center of that cell is an  approximation of the coordinates of the intersection point ${\bm \lambda}$ (see \cite[Section 6]{BMP} and \cite[Section 4]{etal}). 
Of course, such an approximation  is determined up to the chosen discretization of   $\sT $.

Let us stress the fact  that  the results of previous sections \ref{passaggioCella1}, 
\ref{passaggioCella2}, and \ref{passaggioCella2bis} can be used for the construction 
of the accumulator function, as shown in the following algorithm.

\medskip
\noindent
\framebox{
\parbox[l]{0.98\textwidth}
{{\bf The RECOGNITION algorithm}\\
Given, in the image space $\mathbb A^2_{\x}(\reals)$, a Hough regular family $\mathcal F$ 
of irreducible curves of the same degree and a set $\mathbb X = \{p_1,\ldots,p_\nu\}$
of points of interest;
given, in the parameter space $\mathbb A^t_{\bm \Lambda}(\reals)$, 
a region $\mathcal T := [a_1,b_1] \times \cdots \times [a_t,b_t]$  and a discretization, defined by  the initializing point $a=(a_1,\ldots,a_t)\in \sT$  and 
a discretization step $d = (d_1,\ldots,d_t)$,
the algorithm returns a point ${\bm \lambda} \in \mathbb A^t_{\bm \Lambda}(\reals)$.

\begin{description}
\item[1.] For each $p_\iota \in \mathbb X$, let $\mathcal A_\iota$ be the output of the CROSSING AREA 
algorithm applied to the Hough transform $\Gamma_{p_\iota}(\mathcal F)$
w.r.t. to the region $\mathcal T$ and the discretization~$\{a,d\}$.
\item[2.] Compute the multi-matrix $\mathcal A = \sum_{\iota=1}^\nu \mathcal A_\iota$.
\item[3.] In the fixed discretization of $\mathcal T$ find the cell corresponding to 
the unique local maximum of $\mathcal A$; call its center ${\bm \lambda}$ and return it. 
\end{description}
}}

\medskip
We observe  that this approach works for any number $t$ of parameters, while, as far as we know, the Hough transform recognition technique for detection of curves  hardly handles more than three parameters.

In the rest of the section we discuss some illustrative examples,
in which our approach is effectively used to compute 
the accumulator function, which is the core of the recognition algorithm 
based on the Hough transform. In examples \ref{3C} and \ref{EC} below, our outputs are compared with
the results obtained by using well-established pattern recognition techniques for the detection of 
curves in images (see \cite[Sections 6, 7]{BMP} and also \cite[Sections 4, 5]{etal}). Our aim is simply to show that our approach may be successfully used, in a complementary way,  in this context too.
Note that all the computations have been performed on an 
Intel Core 2 Duo processor (at 1.86 GHz), and using the \cocoa 5 implementation 
of the CROSSING AREA algorithm.  We follow  the approach suggested  in Remark \ref{REM2}: in particular, we exploit the bounds $B_1'$ and $B_2'$ from propositions \ref{PropNonPassa2} and \ref{PropPassa2} for degree $\geq 2$, and the bounds $B_1$, $B_2$ from propositions 
\ref{PropNonPassa},  \ref{PropPassa} and Remark \ref{lineare} in the linear case.

All   the families $\sF=\{\sC_{\bm \lambda}\}$ of curves 
in the examples below meet the regularity property  mentioned above  (see  \cite{etal} for details).

\begin{example*}{\bf (Conchoid of Sl\"use)} 
In the affine plane ${\mathbb A}_{(x,y)}^2(\reals)$ consider the  family $\sF=\{\sC_{a,b}\}$ of  rational cubic curves defined by the equation
\begin{equation}\label{CSeq}
\sC_{a,b} : a(x-a)(x^2+y^2) = b^2x^2,
\end{equation}
for some  positive real numbers $a$, $b$. Such a cubic is classically known as 
 {\it conchoid of Sl\"use} of parameters $a$, $b$ (see \cite[p. 130]{atlas}).
The conchoid of Sl\"use is an unbounded rational curve with a 
single singular point  (the origin $O=(0,0)$ is an ordinary double isolated point 
with complex tangent lines of equations $y\pm i\frac{\sqrt{b^2+a^2}}{a}x=0$), 
a single vertical asymptote (the line $x=a$), and a single axis of symmetry 
(the line $y=0$).  Up to the isolated pont $O$, 
the curve $\sC_{a,b}$ lies in the region of the plane
$ \mathbb A^2_{(x,y)}(\mathbb R)$  defined  by  $a < x \le \frac{a^2+b^2}{a}$.

For any point $p=(x(p), y(p))$ of the image plane the Hough transform is  a conic  (an ellipse centered at $\big(\frac{x(p)}{2}, 0\big)$)
in the parameter plane $\mathbb A^2_{(A,B)}(\mathbb R)$ of equation
$$
\Gamma_p(\sF): (x(p)^2+y(p)^2)A^2 + x(p)^2B^2 - x(p)(x(p)^2+y(p)^2)A=0.
$$

In this example we aim to recognize the conchoid of Sl\"use
of parameters $(a,b)=(\frac{1}{4},1)$. To this end, we choose a set $\mathbb X$ 
of $20$ points which lye close to such a curve in the following way. 
We divide the interval $(\frac{1}{4}, \frac{17}{4}]$
in $20$  identical parts: the value of each node represents the $x$-coordinate 
of each point of $\mathbb X$. We then obtain the value of the $y$-coordinates 
by simply solving equation (\ref{CSeq}) in $y$ (with $a=\frac{1}{4}$, $b=1$), 
and picking the rational approximation (with an error of $10^{-1}$) of one of 
its two (symmetric) solutions (this is computed with \cocoa 5, 
using the function {\tt RealRootsApprox}). 
The points of $\mathbb X$ are represented in Figure \ref{conchoidPts}.
\begin{figure}[htb]
\centering
\begin{minipage}[c]{0.48\textwidth}
\includegraphics[width=\textwidth]{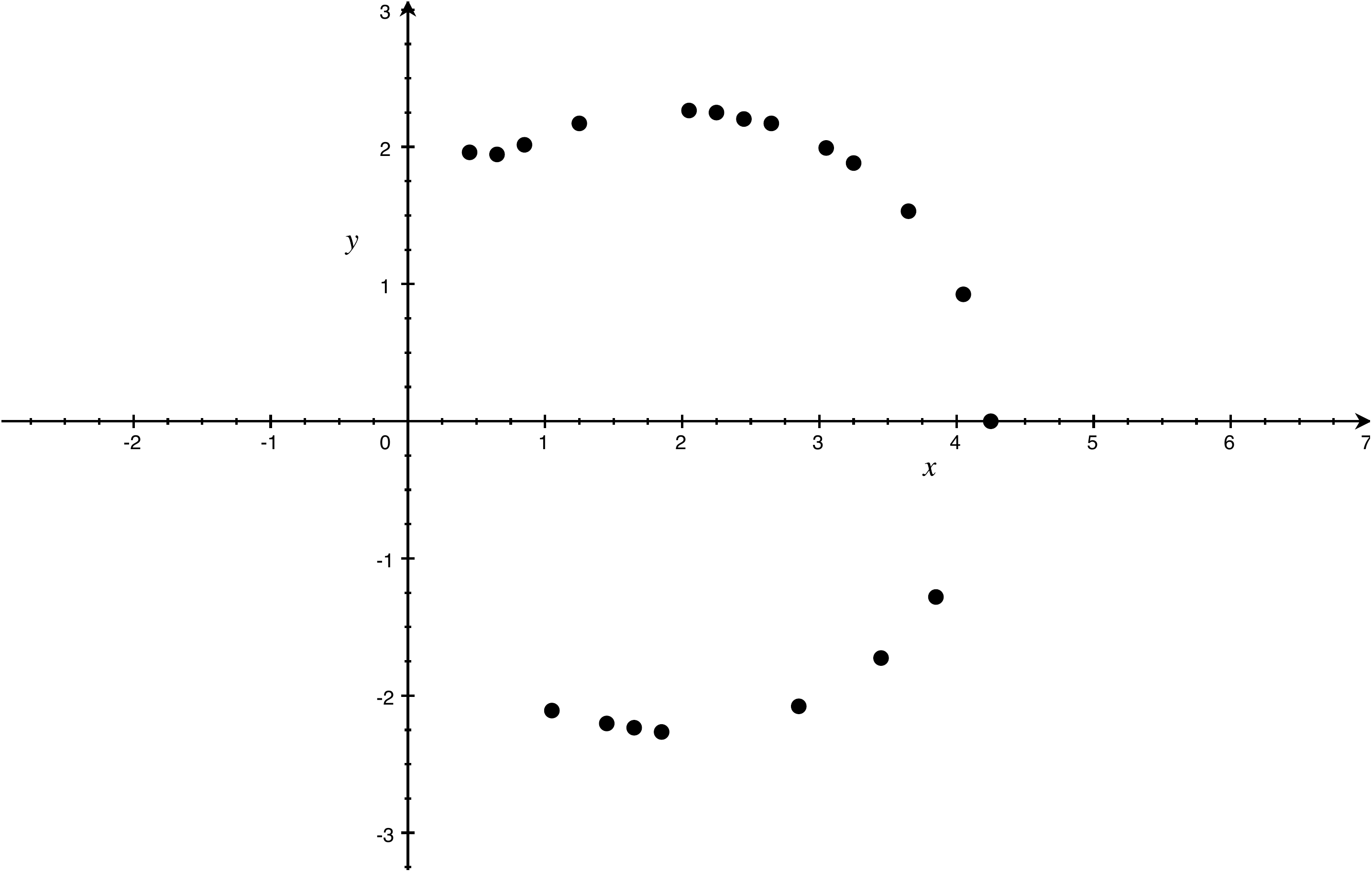}
\caption{\small{The set  $\mathbb X$ of selected points.}}\label{conchoidPts}
\end{minipage}%
\hspace{0.3cm}
\begin{minipage}[c]{0.48\textwidth}
\includegraphics[width=\textwidth]{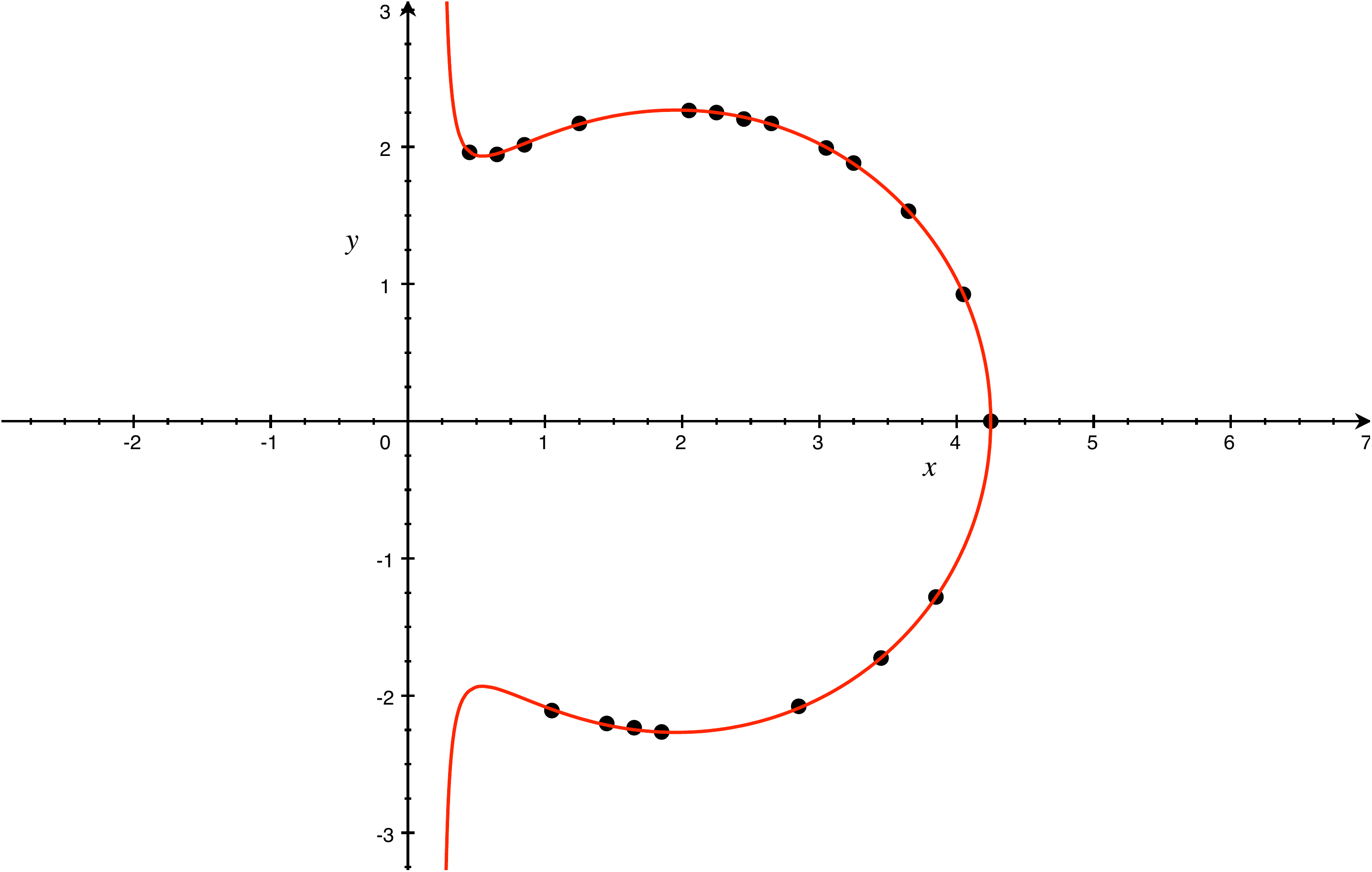}
\caption{\small{The conchoid  $\sC_{\frac{1}{4}, 1}$.}}\label{conchoidCurve}
\end{minipage}
\end{figure}

In the parameter plane, we choose the region $\mathcal T = [0.1,0.5] \times [0.1, 1.1]$
and the discretization step $d=(0.025,0.025)$.
We apply the RECOGNITION algorithm to 
$\mathcal F$, $\mathbb X$, $\mathcal T$, $d$ and we find the following two-dimensional
matrix $\mathcal A$ of size $9 \times 21$:
\begin{eqnarray*}
\small{
\sA = \left ( \begin{array}{ c c c c c c c c c c c c c c c c c c c c c}  
0 & 0 & 0 & 0 & 0 & 0 & 0 & 0 & 0 & 0 & 16 & 18 & 19 & 4 & 2 & 1 & 0 & 0 & 0 & 0 & 0\\
  0 & 0 & 0 & 0 & 0 & 0 & 0 & 0 & 0 & 0 & 0 & 0 & 0 & 15 & 18 & 8 & 1 & 1 & 0 & 0 & 0\\
  0 & 0 & 0 & 0 & 0 & 0 & 0 & 0 & 0 & 0 & 0 & 0 & 0 & 0 & 0 & 9 & 18 & 19 & 1 & 0 & 0\\
  0 & 0 & 0 & 0 & 0 & 0 & 0 & 0 & 0 & 0 & 0 & 0 & 0 & 0 & 0 & 0 & 0 & 0 & \bf{20} & 0 & 0\\
  0 & 0 & 0 & 0 & 0 & 0 & 0 & 0 & 0 & 0 & 0 & 0 & 0 & 0 & 0 & 0 & 0 & 1 & 0 & 19 & 16\\
  0 & 0 & 0 & 0 & 0 & 0 & 0 & 0 & 0 & 0 & 0 & 0 & 0 & 0 & 1 & 1 & 1 & 0 & 0 & 1 & 4\\
  0 & 0 & 0 & 0 & 0 & 0 & 0 & 0 & 0 & 0 & 1 & 1 & 1 & 1 & 0 & 0 & 0 & 0 & 1 & 0 & 1\\
  1 & 1 & 1 & 1 & 1 & 1 & 1 & 1 & 0 & 0 & 0 & 0 & 0 & 0 & 0 & 0 & 0 & 1 & 0 & 0 & 1\\
  0 & 0 & 0 & 0 & 0 & 0 & 0 & 0 & 0 & 0 & 0 & 0 & 0 & 0 & 0 & 1 & 1 & 0 & 0 & 0 & 1
 \end{array} \right )
 }
\end{eqnarray*}
The maximum entry of $\sA$ is $\bf{20}$ (which is exactly the cardinality of $\mathbb X$), and it 
corresponds to the cell centered at $(A,B)= (\frac{1}{4},1)$, that is, exactly to 
the point of the parameter plane which identifies the curve~$\sC_{\frac{1}{4}, 1}$ 
we started from.
\end{example*}

The  examples below  consist of an application
of the recognition algorithm to detect profiles of interest in medical imaging. They comes from  X-ray Computerized Tomography images  studied in \cite{etal}: see in particular Figure 3 (detection of a lumbar vertebra profile) and Figure 4 (detection  of the canal spynal) of that paper.

\begin{example*}\label{3C}{\bf (Curve with $3$ convexities)}
In the affine plane ${\mathbb A}_{(x,y)}^2(\reals)$ consider the family 
$\sF=\{\mathcal{C}_{a,b}\}$ of sextic curves
  \begin{equation*}\label{3convex}
\mathcal{C}_{a,b}: (x^2+y^2)^3 = (a(x^2+y^2) - b(x^3-3xy^2))^2,
\end{equation*}
for positive real numbers $a$, $b$, with $b<1$. 
Such curves, classically known as {\em curve with $3$ convexities}
 (see \cite[p. 183]{atlas}),
are bounded and contained in the circular crown of radii $\frac{a}{1+b}$
and $\frac{a}{1-b}$.
The origin $O=(0,0)$ is an isolated point with multiplicity $4$, 
with complex conjugates  tangent lines of equation 
$(x^2+y^2)^2=0$. Further, $O$ has some special complexity:
a  more detailed local study of the curve at $O$ shows that in fact $\sC_{a,b}$ 
is rational.

For any point $p=(x(p), y(p))$ of the image plane the Hough transform is the  
degenerate conic (i.e., the union of two parallel lines) in the parameter plane $
{\mathbb A}_{(A,B)}^2(\reals)$ of equation
$$
\Gamma_p(\sF): (x(p)^2+y(p)^2)^3 = (A(x(p)^2+y(p)^2) - B(x(p)^3-3x(p)y(p)^2))^2.
$$

In the image plane ${\mathbb A}_{(x,y)}^2(\reals)$ we consider $1170$ points $\mathbb Y$ represented in Figure \ref{3convexBottomPoints} (the data set $\mathbb Y$  has been taken  from \cite{etal}, see Section 5 and Figure 4, after applying a standard edge detection algorithm). We zoom in the image to consider the portion
of interest, and extract from $\mathbb Y$ the set of points lying in the box $[-1.5,1.5] \times [-1.5,1.5]$.
Such a set  is denoted by $\mathbb X$ and consists of $320$ points,
represented in Figure \ref{3convexBottomCurve}.
Taking into account the variance of the curve, we choose in the parameter plane 
the region $\mathcal T =[0.35,0.9] \times [0.175,0.5] $
and the discretization step $d=(0.015, 0.015)$.
Applying the RECOGNITION algorithm to 
$\mathcal F$, $\mathbb X$, $\mathcal T$, $d$ we get a two-dimensional
matrix $\mathcal A$ of size $19 \times 12$.
The (unique) maximum entry of $\sA$ 
corresponds to the cell centered at  $(A,B) = (0.53,0.445)$,
which yields the red curve represented in Figure \ref{3convexBottomCurve}.

\begin{figure}[htb]
\centering
\begin{minipage}[c]{0.51\textwidth}
\includegraphics[width=\textwidth]{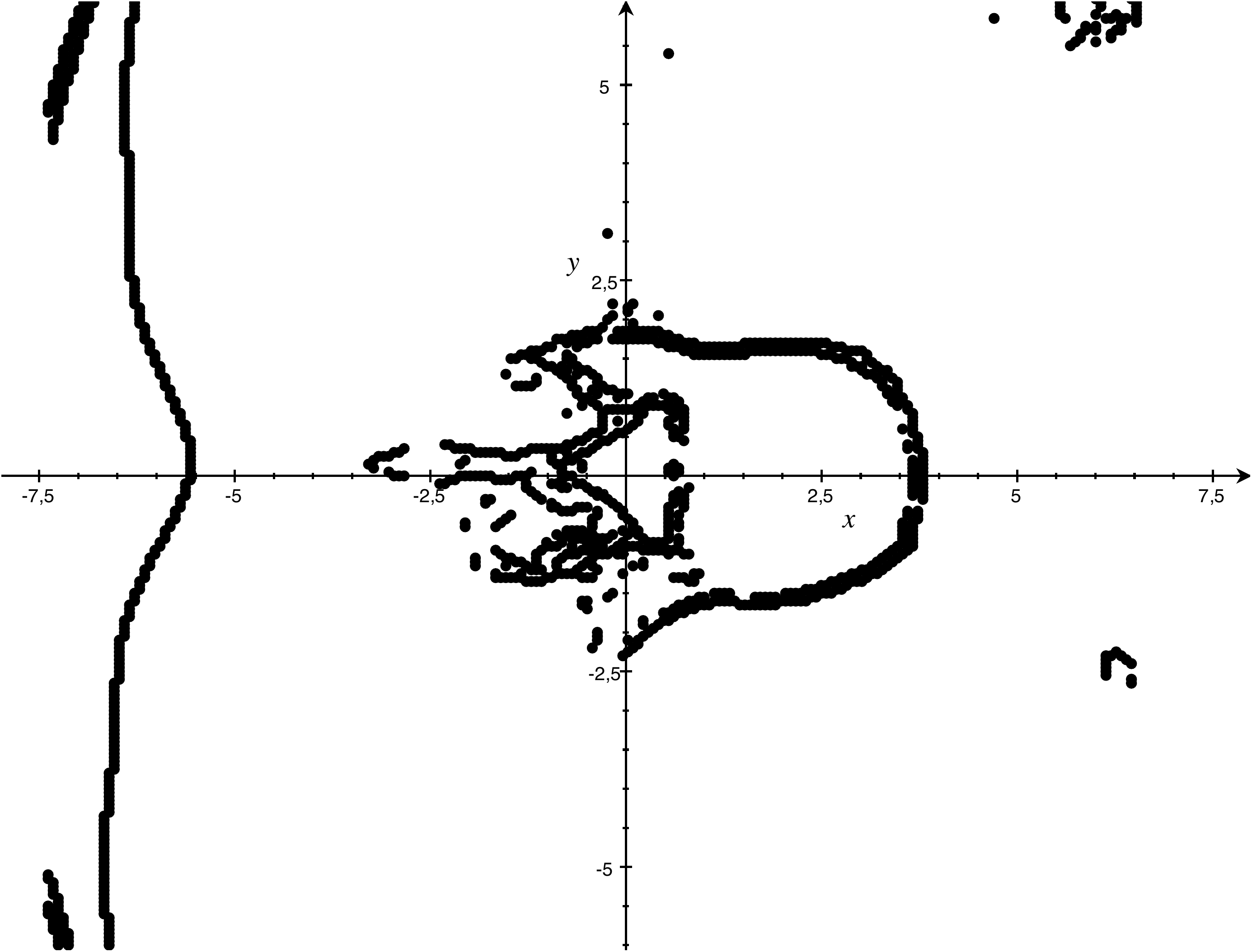}
\caption{\small{The points of the set $\mathbb Y$.}}\label{3convexBottomPoints}
\end{minipage}%
\hspace{0.29cm}
\begin{minipage}[c]{0.46\textwidth}
\includegraphics[width=\textwidth]{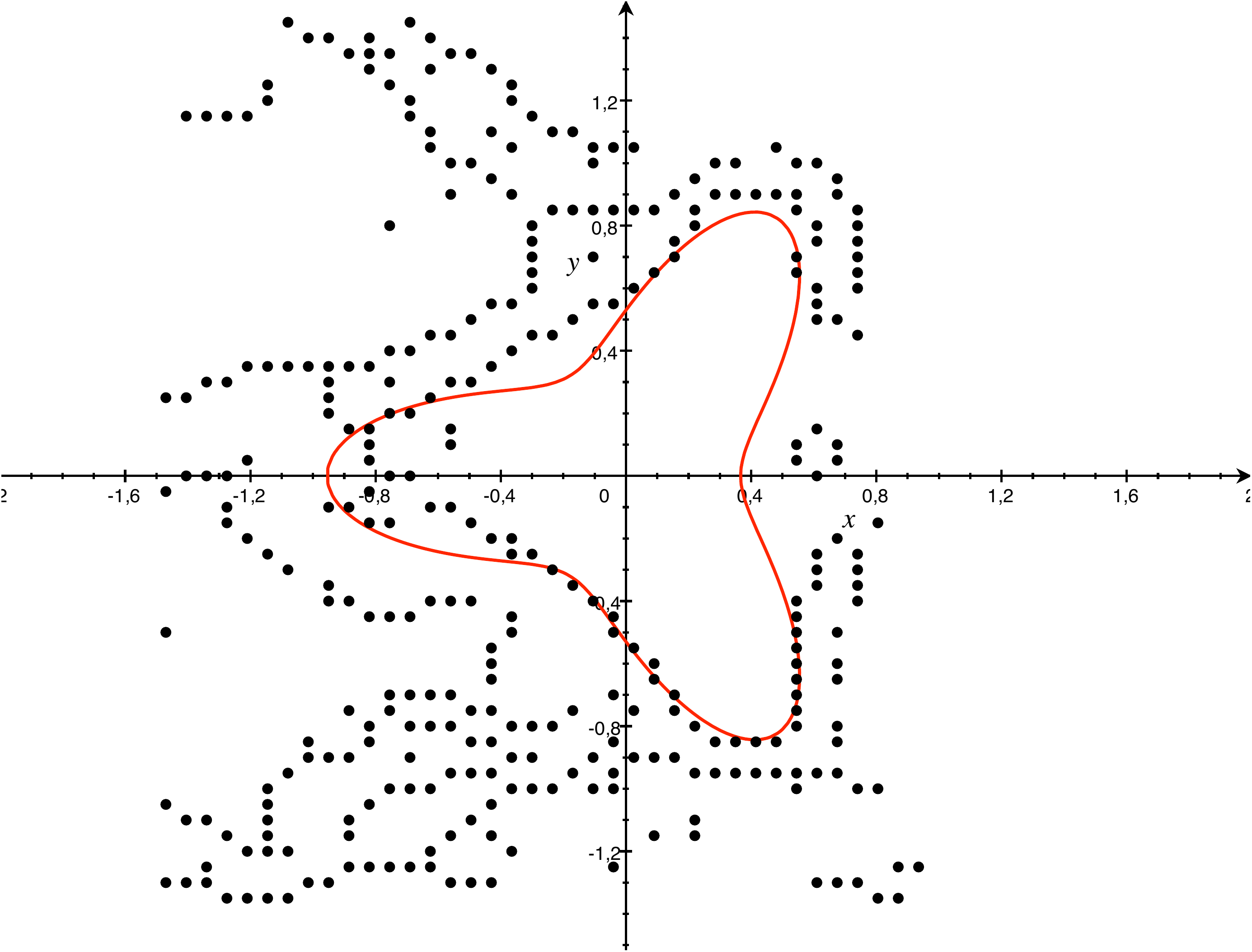}
\caption{\small{The points of the set $\mathbb X$ and the curve with $3$ convexities
$\sC_{0.53,0.445}$: focus on the portion of interest.}}\label{3convexBottomCurve}
\end{minipage}
\end{figure}

Consider now the family $\sF=\{\mathcal{C}_{a,b,m}\}$ 
of sextic curves
  \begin{equation*}\label{diludendo}
\mathcal{C}_{a,b,m}: (mx^2+y^2)^3 = (a(mx^2+y^2) - b(x^3-3xy^2))^2,
\end{equation*}
 for  real positive parameters $a$, $b$,  $m$, with $b<1$. 
 This is a slight variant of  the curve with $3$ convexities
 which corresponds to the case $m=1$. 
 
 For any point $p=(x(p), y(p))$ of the image plane the Hough transform is the  
quartic  surface  in the parameter space ${\mathbb A}_{(A,B,M)}^3(\reals)$ of equation
$$
\Gamma_p(\sF): (Mx(p)^2+y(p)^2)^3 = (A(Mx(p)^2+y(p)^2) - B(x(p)^3-3x(p)y(p)^2))^2.
$$

In the image plane ${\mathbb A}_{(x,y)}^2(\reals)$  we consider  $2433$ points  $\mathbb Y$ represented in Figure \ref{3convexTopPoints}  
(the data set  $\mathbb Y$ has  been taken from \cite[Figure 4, Top case]{etal}, after an edge detection processing). 
We zoom in the image to consider the portion
of interest, and extract from $\mathbb Y$  the set of points lying in the box $[-1.5,1.5] \times [-1.5,1.5]$.
Such a set is denoted by $\mathbb X$ and consists of $132$ points,
represented in Figure \ref{3convexTopCurve}.
Taking into account the variance of the curve, we choose in the parameter plane 
the region $\mathcal T =[0.7,1] \times [0,  0.18] \times [0.9,1.1] $
and the discretization step $d=(0.02, 0.02, 0.02)$.
Applying the RECOGNITION algorithm to 
$\mathcal F$, $\mathbb X$, $\mathcal T$, $d$ we get  a  three-dimensional
matrix $\mathcal A$ of size $9 \times 6 \times 6$.
The (unique) maximum entry of $\sA$  
corresponds to the cell centered at $(A,B,M) = (0.82, 0.04, 1.1)$
which yields the red curve represented in Figure \ref{3convexTopCurve}.

\begin{figure}[htb]
\centering
\begin{minipage}[c]{0.50\textwidth}
\includegraphics[width=\textwidth]{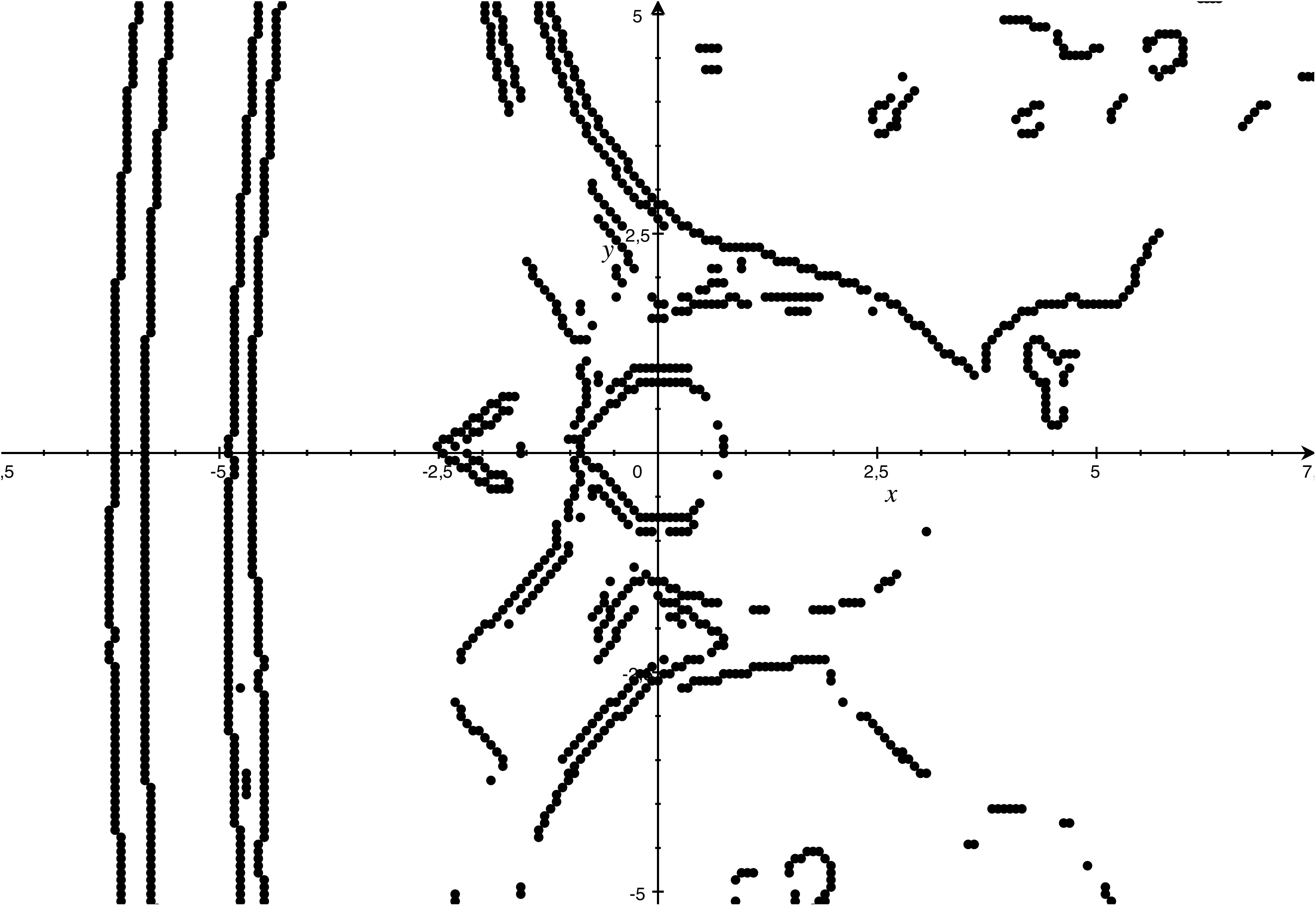}
\caption{\small{The points of the set $\mathbb Y$.}}
\label{3convexTopPoints}
\end{minipage}%
\hspace{0.29cm}
\begin{minipage}[c]{0.47\textwidth}
\includegraphics[width=\textwidth]{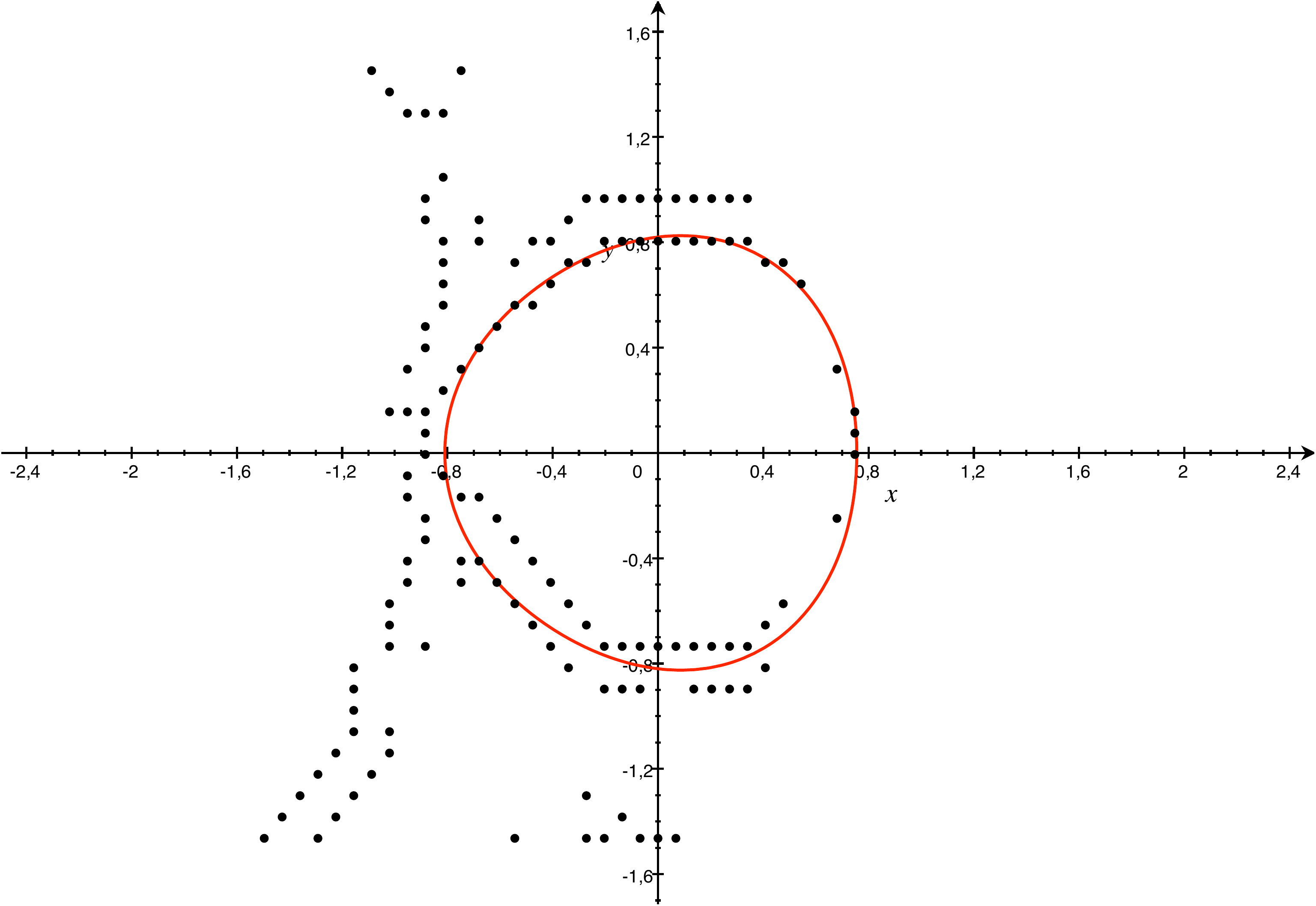}
\caption{\small{The points of the set $\mathbb X$ and the curve with $3$ convexities
$\sC_{0.82, 0.04, 1.1}$: focus on the portion of interest.}}\label{3convexTopCurve}
\end{minipage}
\end{figure}
\end{example*}

\begin{example*}\label{EC}{\bf (Elliptic curves)} 
In the affine plane ${\mathbb A}_{(x,y)}^2(\reals)$ consider the family $\sF=\{ \mathcal{C}_{a,b,m,n}\}$ 
of cubic curves of equation
  $$
\mathcal{C}_{a,b,m,n}: x^2=-my^3 + ny^2-ay + b, 
$$
for real parameters $a$, $b$, $m$, $n$ with $m$ positive.
The general curve of the family is non-singular, so that it is an elliptic curve. 
E.g., for $m=1$ and $n=0$,  we get the so called {\em Weierstrass equation}. 
For special values of the parameters $a$, $b$, $m$, $n$, the cubic  has either 
a nodal or a cuspidal double point (as clearly  it happens e.g. for $a=b=n=0$). 

For any point $p=(x(p), y(p))$ of the image plane the Hough transform is the 
hyperplane  in the parameter space ${\mathbb A}_{(A,B,M,N)}^4(\reals)$ of equation
$$
\Gamma_p(\sF): My(p)^3-Ny(p)^2+Ay(p) -B+x(p)^2=0.
$$

In the image plane ${\mathbb A}_{(x,y)}^2(\reals)$ we consider $1084$ points $\mathbb Y$ represented in Figure \ref{ellitticaPoints}
(the data set $\mathbb Y$ has been taken from \cite{etal}, see Figure 3 (upper panels), after an edge detection processing).
We zoom in the image to consider the portion
of interest, and extract from $\mathbb Y$ the  set of points lying in the box $[-3,3] \times [-3, 3]$.
Such a set  is denoted by $\mathbb X$ and consists of $636$ points,
represented in Figure \ref{ellitticaCurve}.
After looking at the variance of the curve, we choose the region  in the parameter space
$\mathcal T =[-1.02,0.206] \times [1.96,2.89]  \times [0.8,1.2] \times [-0.2,0.2]$
and the discretization step $d=(0.1, 0.1, 0.1, 0.1)$.
We apply the RECOGNITION algorithm to 
$\mathcal F$, $\mathbb X$, $\mathcal T$, $d$ and find a four-dimensional
matrix $\mathcal A$ of size $7 \times 6 \times 3 \times 3$ which exhibits 
a unique  maximum. Its value 
corresponds to the cell centered at $(A,B,M,N) =(-0.42, 2.76, 0.8, 0)$
which yields the red curve represented in Figure \ref{ellitticaCurve}.

\begin{figure}[htb]
\centering
\begin{minipage}[c]{0.50\textwidth}
\includegraphics[width=\textwidth]{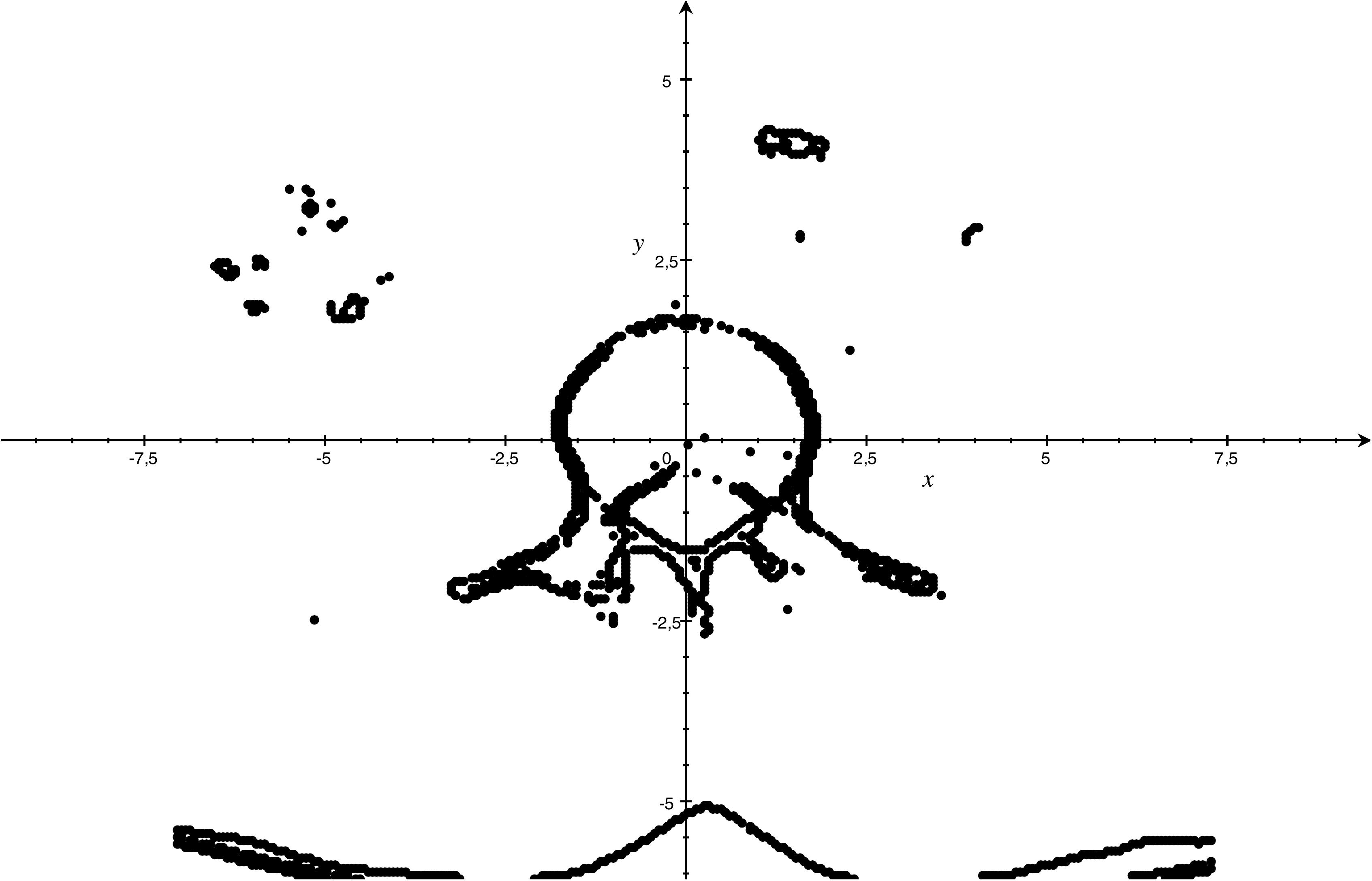}
\caption{\small{The points of the sets $\mathbb Y$.} }
\label{ellitticaPoints}
\end{minipage}%
\hspace{0.29cm}
\begin{minipage}[c]{0.47\textwidth}
\includegraphics[width=\textwidth]{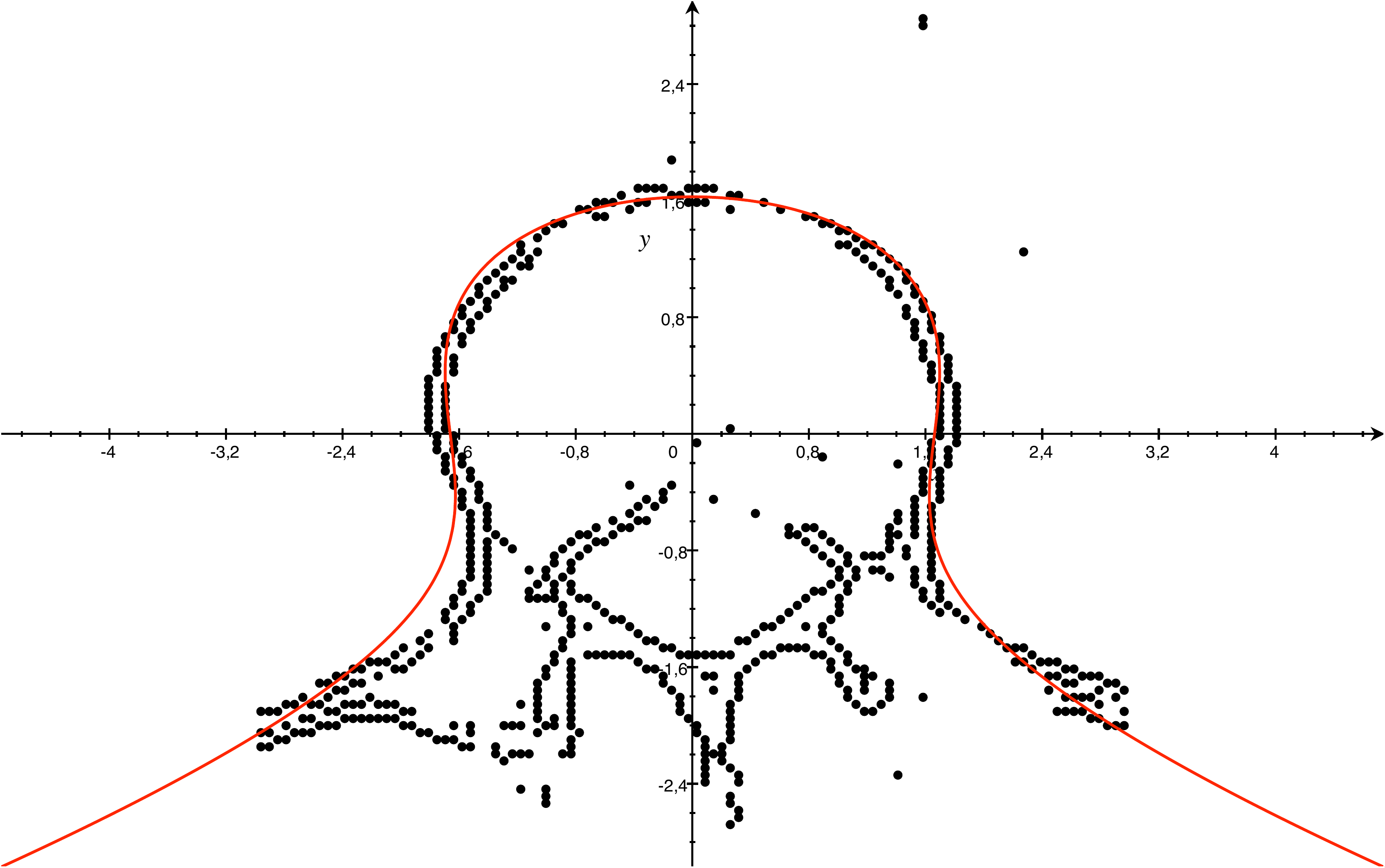}
\caption{\small{The points of the set $\mathbb X$ and the elliptic curve 
$\sC_{-0.42, 2.76, 0.8, 0}$: focus on the portion of interest.}}\label{ellitticaCurve}
\end{minipage}
\end{figure}
\end{example*}

In closing,  we propose a comparison attempt of  the outputs of our algorithm   and  the  pattern recognition techniques used in  \cite{etal},  where curves from the families $\mathcal F$ studied in examples \ref{3C} and \ref{EC} have been used  to detect vertebrae profiles. The black curves in figures \ref{3convexBottomCompCurves} and \ref{ellitticCompCurves} below  are those detected  in Figures 4 (Bottom case) and Figure 3 (panel  (d))  of \cite{etal}, respectively.

\begin{figure}[htb]
\centering
\begin{minipage}[c]{0.4\textwidth}
\includegraphics[width=\textwidth]{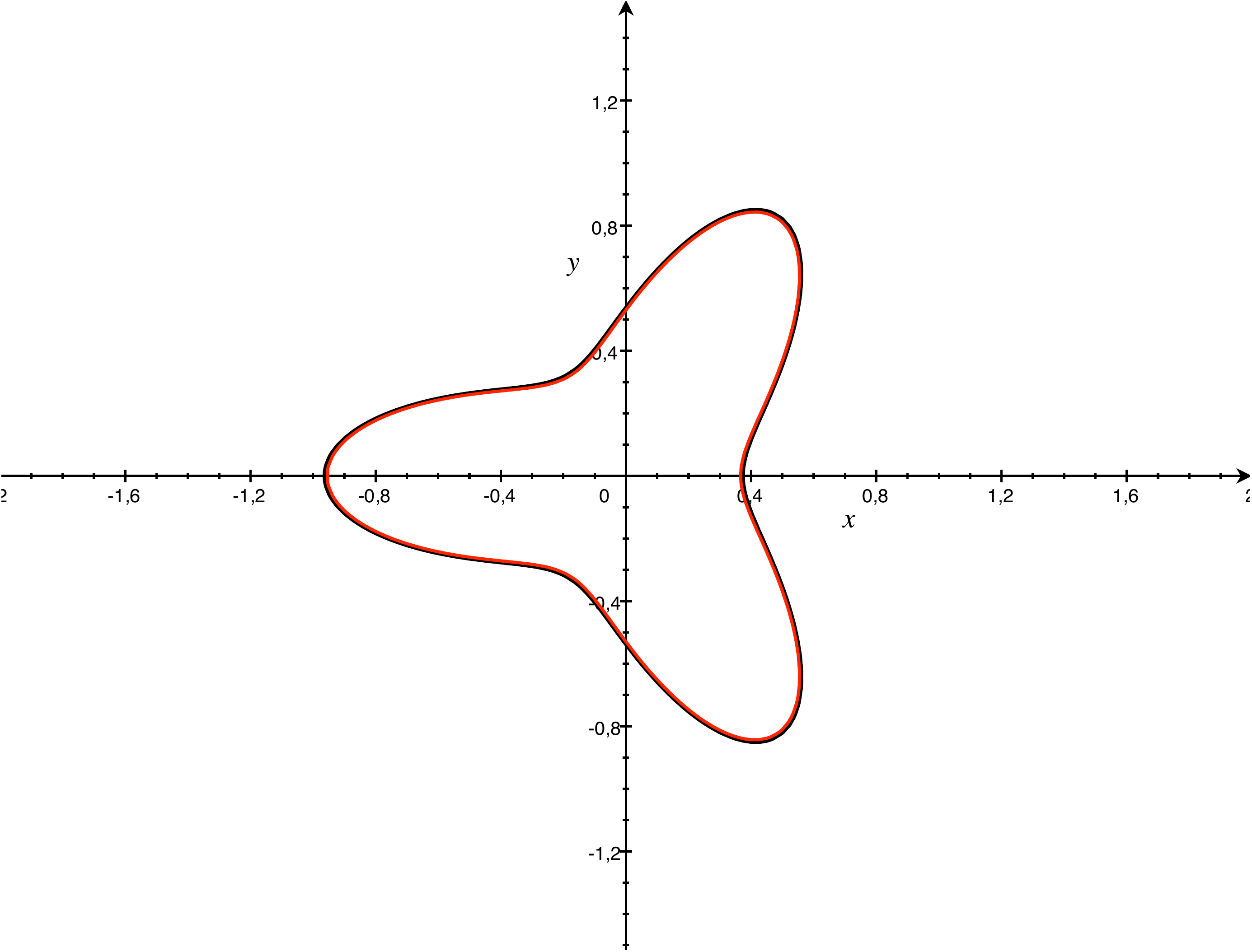}
\caption{\small{The curves with $3$ convexities: $\sC_{0.53,0.445}$ (red,  Figure \ref{3convexBottomCurve}) and $\sC_{0.54, 0.44}$ (black).}}\label{3convexBottomCompCurves}
\end{minipage}%
\hspace{0.29cm}
\begin{minipage}[c]{0.47\textwidth}
\includegraphics[width=\textwidth]{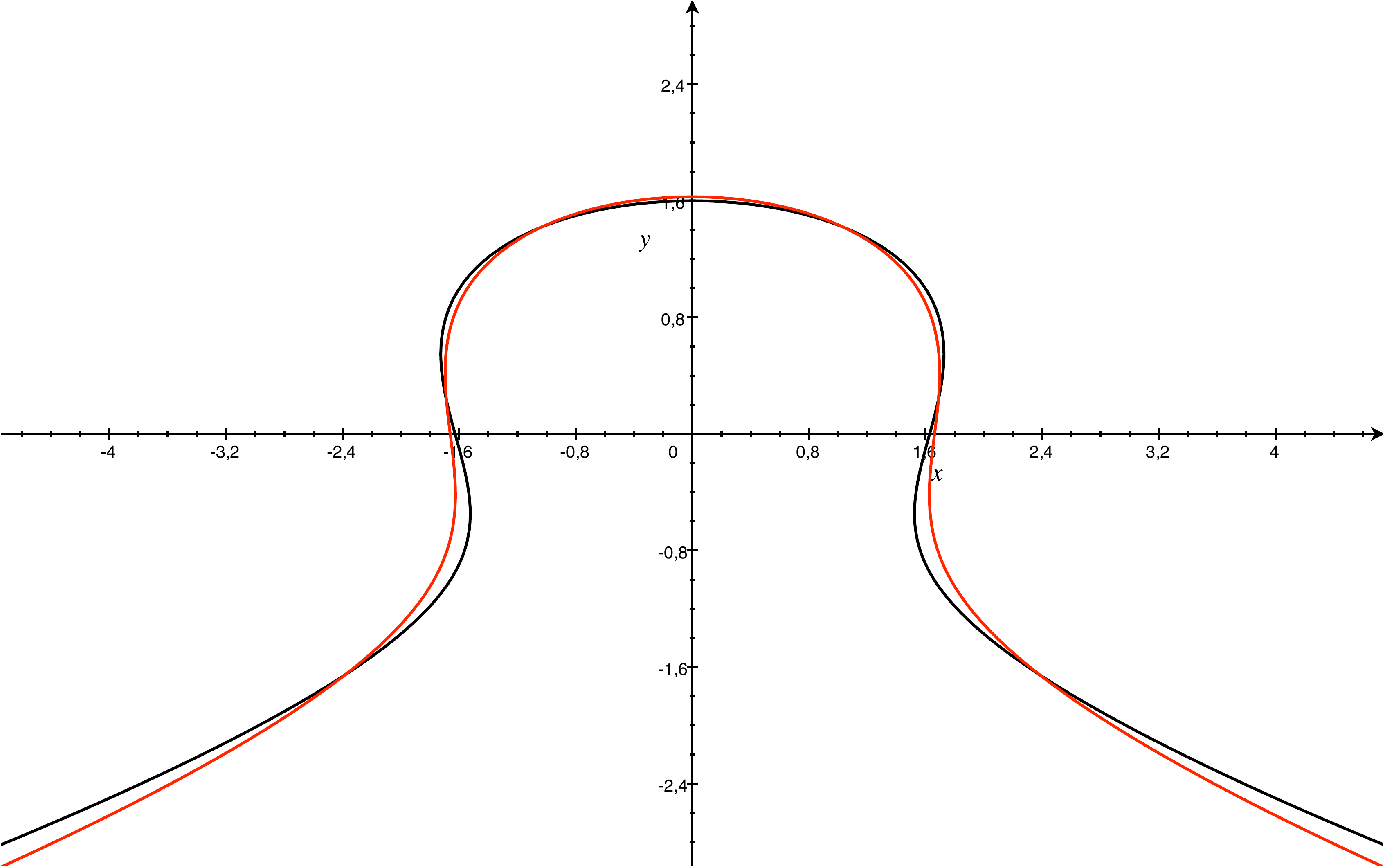}
\caption{\small{The elliptic curves: $\sC_{-0.42, 2.76, 0.8, 0}$ (red, Figure \ref{ellitticaCurve}) and $\sC_{-0.9, 2.65, 1, 0}$ (black).}}\label{ellitticCompCurves}
\end{minipage}
\end{figure}

\small{

}

\bigskip
\bigskip

\noindent
 M. Torrente,
 Dipartimento di Matematica,
Universit\`a di Genova,
Via Dodecaneso 35,
I-16146 Genova, Italy. e-mail {\tt torrente@dima.unige.it}

\smallskip

\noindent M.C. Beltrametti,
Dipartimento di Matematica,
Universit\`a di Genova,
Via Dodecaneso 35,
I-16146 Genova, Italy. e-mail {\tt beltrame@dima.unige.it}